\documentclass[final,1p,times]{elsarticle}
\usepackage{amsmath,amssymb,amsfonts,amsthm} 
\usepackage{hyperref}                 
\newtheorem{theorem}{Theorem}[section]
\newtheorem{proposition}[theorem]{Proposition}
\newtheorem{corollary}[theorem]{Corollary}
\newtheorem{lemma}[theorem]{Lemma}
\newtheorem{remark}[theorem]{Remark}

\newtheorem{example}[theorem]{Example}
\numberwithin{equation}{section}

\journal{Journal of Number Theory}

\begin{document}

\begin{frontmatter}

\title{A Gauss-Kuzmin-type problem for a family of continued fraction expansions}

\author[anmb]{Dan Lascu}
\ead{lascudan@gmail.com}
\address[anmb]{Mircea cel Batran Naval Academy, 1 Fulgerului, 900218 Constanta,Romania}

\begin{abstract}
In this paper we study in detail a family of continued fraction expansions of any number in the unit closed interval $[0,1]$ whose digits are differences of consecutive non-positive integer powers of an integer $m \geq 2$. For the transformation which generates this expansion and its invariant measure, the Perron-Frobenius operator is given and studied. For this expansion, we apply the method of random systems with complete connections by Iosifescu and obtained the solution of its Gauss-Kuzmin type problem.
\end{abstract}
\end{frontmatter}

\sloppy

\section{Introduction}

The purpose of this paper is to prove a Gauss-Kuzmin type problem for non-regular continued fraction expansions introduced by Chan \cite{r5}. In order to solve the problem, we apply the random systems with complete connections by Iosifescu \cite{r10}. First we outline the historical framework of this problem. Then, in Section 1.2, we present the current framework. The main theorem will be shown in Section 1.3. In this subsection we will also give a detailed outline of the paper.

\subsection{Gauss' Problem}

One of the first and still one of the most important results in the metrical theory of continued fractions is so-called Gauss-Kuzmin theorem. Write $x \in [0,1)$ as a regular continued fraction 
\[
x = \displaystyle \frac{1}{a_1+\displaystyle \frac{1}{a_2+\displaystyle \frac{1}{a_3+ \ddots}}} :=[a_1, a_2, a_3, \ldots],
\]
where $a_n \in \mathbb{N}_+ : = \left\{1, 2, 3, \ldots\right\}$.
The metrical theory of continued fractions started on 25th October 1800, with a note by Gauss in his mathematical diary. Gauss wrote that (in modern notation) 
\[
\lim_{n \rightarrow \infty} \lambda \left(\tau^n \leq x\right) = \frac{\log(1+x)}{\log2}, \ x \in I:=[0,1].
\]
Here $\lambda$ is Lebesgue measure and the map $\tau : [0, 1) \rightarrow [0, 1)$, the so-called \textit{regular continued fraction} (or \textit{Gauss}) transformation, is defined by 
\[
\tau(x) := \frac{1}{x}-\left\lfloor \frac{1}{x} \right\rfloor, \quad x \neq 0; \ \tau(0):=0,
\]
where $\left\lfloor \cdot \right\rfloor$ denotes the \textit{floor} (or \textit{entire}) function. Gauss' proof (if any) has never been found. A little more than 11 years later, in a letter dated 30 January 1812, Gauss asked Laplace to estimate the error
\[
e_n(x) := \lambda \left(\tau^{-n}[0, x]\right) - \frac{\log(1+x)}{\log2}, \quad n \geq 1, \ x\in I.
\]
This has been called \textit{Gauss' Problem}. It received a first solution more than a century later, when R.O. Kuzmin (see \cite{r17}) showed in 1928 that $e_n(x) = \mathcal{O}(q^{\sqrt{n}})$ as $n \rightarrow \infty$, uniformly in $x$ with some (unspecified) $0 < q < 1$. One year later, using a different method, Paul L\'evy (see \cite{r19}) improved Kuzmin's result by showing that $\left|e_n(x)\right| \leq q^n$, $n \in \mathbb{N}_+$, $x \in I$, with $q = 3.5 - 2\sqrt{2} = 0.67157...$. The Gauss-Kuzmin-L\'evy theorem is the first basic result in the rich metrical theory of continued fractions. 

\subsection{A non-regular continued fraction expansion}

In this paper, we consider a generalization of the Gauss transformation and prove an analogous result.

In \cite{r5}, Chan shows that any $x \in \left[0, 1\right)$ can be written in the form
\begin{equation}
x = \frac{m^{-a_1(x)}}{\displaystyle 1+\frac{(m-1)m^{-a_2(x)}}{\displaystyle 1 + \frac{(m-1)m^{-a_3(x)}}{\displaystyle 1 + \ddots} }}:=[a_1(x), a_2(x), a_3(x), \ldots]_m, \label{3.1}
\end{equation}
where $m \in \mathbb{N}_+$, $m \geq 2$ and $a_n(x)$'s are non-negative integers.

For any $m \in \mathbb{N}_+$ with $m \geq 2$, define the transformation $\tau_m$ on $I$ by
\begin{equation}
\tau_m (x) = \left\{\begin{array}{lll} 
\displaystyle \frac{m^{\left\{\frac{\log x^{-1}}{\log m}\right\}}-1}{m-1}, & \hbox{if} & x \ne 0 \\ 
\\
0, & \hbox{if} & x = 0,
\end{array} \right. \label{3.5}
\end{equation}
where $\left\{\cdot\right\}$ stands for fractionary part. It is easy to see that $\tau_m$ maps the set $\Omega$ of irrationals in $I$ into itself. For any $x \in (0,1)$ put
\begin{equation}
a_n=a_n(x) = a_1\left(\tau_m^{n-1}(x)\right), \quad n \in {\mathbb{N}}_+, \label{1.3}
\end{equation}
with $\tau_m^0 (x) = x$ and
\begin{equation}
a_1=a_1(x) = \left\{\begin{array}{lll} 
\lfloor\log x^{-1} / \log m\rfloor, & \hbox{if} & x \neq 0 \\ 
\infty, & \hbox{if} & x = 0.
\end{array} \right. \label{1.4}
\end{equation}
Transformation $\tau_m$ which generates the continued fraction expansion (\ref{3.1}) is ergodic with respect to an invariant probability measure, $\gamma_m$, where 
\[
\gamma_m (A) = k_m \int_{A} \frac{dx}{((m-1)x+1)((m-1)x+m)}, \quad A \in {\mathcal{B}}_I,
\]
with $k_m = \frac{(m-1)^2}{\log \left(m^2/(2m-1)\right)}$ and $\mathcal{B}_I$ is the $\sigma$-algebra of Borel subsets of $I$ (which, by definition, is the smallest $\sigma$-algebra containing intervals).

The ergodicity of $\tau_m$ plays a key role in the study of the asymptotic growth rate of the random Fibonacci type sequences $\left\{f_n\right\}$ defined by $f_{-1} = 0$, $f_0 = 1$, $c_0 = 0$ and
\begin{equation}
f_{n}=m^{c_n} f_{n-1}+(m-1)m^{c_{n-1}}f_{n-2}, \label{1.3'}
\end{equation}
where $c_n$, $n \geq 1$, are the digits from (\ref{3.1}). As is known, the Fibonacci sequence is defined using the linear recurrence relation
\[
F_{n+1} = F_n + F_{n-1},  \ n \in \mathbb{N}_+, \mbox{with} \ F_0 = F_1 = 1,
\]
and Binet's formula is 
\[
F_n = \displaystyle\frac{1}{\sqrt{5}} \left(\displaystyle\frac{1+\sqrt{5}}{2}\right)^{n+1} - \displaystyle\frac{1}{\sqrt{5}} \left(\displaystyle\frac{1-\sqrt{5}}{2}\right)^{n+1}, \ n \in \mathbb{N}. 
\]
It is known that using Binet's formula we can compute the asymptotic growth rate of the Fibonacci sequence $\left\{F_n\right\}$, which is given by
\[
\lim_{n \rightarrow \infty} \frac{1}{n} \log F_n = \log \left(\frac{1+\sqrt{5}}{2}\right) = 0.4812\ldots
\]
In the case of random Fibonacci type sequences, defined by (with fixed $f_1$ and $f_2$) 
\[
f_n = \alpha(n) f_{n-1} + \beta(n) f_{n-2},
\]
where $\alpha(n)$ and $\beta(n)$ are random coefficients, the quest for the asymptotic growth rate is more difficult. Recently, Viswanath (see \cite{Viswanath-2000}) proved that the asymptotic growth rate of the random Fibonacci sequences defined by $f_1 = f_2 = 1$ and 
\[
f_n = \pm f_{n-1} \pm f_{n-2},
\]
where the signs are chosen independently and with equal probabilities, is given by
\[
\lim_{n \rightarrow \infty} \frac{1}{n} \log f_n = \log (1.13198824\ldots) = 0.12397559\ldots
\]
with probability $1$. But Viswanath's method is not the only way through. So, Chan proved in \cite{r5} that for almost all $x$ with respect to the Lebesgue measure, the asymptotic growth rate of $\{f_n\}$ from (\ref{1.3'}) is given by
\begin{eqnarray*}
\lim_{n \rightarrow \infty} \frac{1}{n} \log f_n &=& k_m \int^{1}_{0} \frac{\log(1/x)}{((m-1)x + 1)((m-1)x + m)}\mbox{d}x \\
&\leq& k_m \frac{3m-1}{2m(2m-1)}. \qquad \qquad \qquad \qquad \qquad \qquad \qquad \ \ \ \ \ \ \ \ \ \ \hfill \Box 
\end{eqnarray*}

\subsection{Main theorem}

We show our main theorem in this subsection. For this purpose let $\mu$ be a non-atomic probability measure on $\mathcal{B}_I$ and define 
\begin{eqnarray*}
F_n (x) &=& \mu (\tau_m^n < x), \ x \in I, \ n \in \mathbb{N}, \\
F(x) &=& \displaystyle \lim_{n \rightarrow \infty}F_n(x), \ x \in I,
\end{eqnarray*}
with $F_0 (x) = \mu ([0,x))$. 

Then the following holds.
\begin{theorem} $\mathrm{(A \ Gauss-Kuzmin-type \ theorem)}$ \label{G-K}
If $\mu$ has a Riemann-integrable density, then
\begin{equation}
F(x) = \frac{k_m}{(m-1)^2}\log \frac{m((m-1)x+1)}{(m-1)x+m}, \quad x \in I, \label{86}
\end{equation}
where $k_m = \displaystyle \frac{(m-1)^2}{\log \left(m^2/(2m-1)\right)}$.

If the density of $\mu$ is a Lipschitz function, then there exist two positive constants $q < 1$ and $k$ such that for all $x \in I$ and $n \in \mathbb{N}_+$ we have
\begin{equation}
\mu \left(\tau^n_m < x\right) = \frac{k_m}{(m-1)^2}(1 + \theta q^n) \log \frac{m((m-1)x+1)}{(m-1)x+m}, \label{87}
\end{equation}
where $\theta$ is a certain constant determined by $\mu,n,x$ such that $|\theta|\leq k$.
\end{theorem}

The paper is organised as follows. In Section 2 we give the basic metric properties of the continued fraction expansion in (\ref{3.1}). Hence, we give a Legendre-type result and the Brod\'en-Borel-L\'evy formula used to determine the probability structure of $(a_n)_{n \in \mathbb{N}_+}$ under $\lambda$. In Section 2.4, we find the invariant measure of $\tau_m$. The proof of this result is given in a different manner from that described by Chan in \cite{r5}. In Section 3 we consider the so-called \textit{natural extension} $\overline{\tau}_m$ (see \cite{r20}), define extended incomplete quotients $\overline{a}_l$, $l \in \mathbb{Z}$ and we generalize some results presented in Section 2. In Section 4, we derive the associated Perron-Frobenius operator under different probability measures on ${\mathcal{B}}_I$. We study the Perron-Frobenius operator of $\tau_m$ under the invariant measure $\gamma_m$ induced by the limit distribution function, we derive the asymptotic behaviour of this operator and we restrict the Perron-Frobenius operator to the linear space of all complex-valued functions of bounded variation and to the space of all bounded measurable complex-valued functions. Section 5 is divided into three parts. The first subsection has as purpose defining the notion of random system with complete connections. In the second subsection we set up the necessary machinery to prove the main theorem whose proof is contained in the last subsection. To determine where $\mu (\tau_m^n < x)$ tends as $n \rightarrow \infty$ and give the rate of this convergence, we use the ergodic behaviour of the random system with complete connections associated with this expansion. For a more detailed study of the theory and applications of dependence with complete connections to the metrical problems and other interesting aspects of number theory we refer the reader to \cite{r10,r13,r14,Sebe-2001,Sebe-2002,r23} and others. 

\section{Metric properties of the continued fraction expansions in (\ref{3.1})}

Roughly speaking, the metrical theory of continued fraction expansions is about properties of the sequence $(a_n)_{n \in \mathbb{N}}$ and related sequences (see section 3.2). 
The main purpose of this section is to determine the probability structure of $(a_n)_{n \in \mathbb{N}_+}$ under the Lebesgue measure $\lambda$. Before that, we shortly present the metrical theory of these continued fraction expansions. Another important result is the \textit{Legendre's theorem-type} (see, e.g., \cite{BK-1990, r11, Hincin-1964}) which is one of the main reasons for studying continued fractions, because it tells us that good approximations of irrational numbers by rational numbers are given by continued fraction convergents. 

\subsection{Some elementary properties of the continued fraction expansion in (\ref{3.1})}

Here, we want to prove the convergence of expansion of the type of (\ref{3.1}). First, note that in the rational case, the continued fraction expansion (\ref{3.1}) is finite, unlike the irrational case, when we have an infinite number of non-negative digits. 

Define $[a_1, a_2, \ldots, a_n]_m$ the convergent of $\omega \in \Omega$ by truncating the expansion on the right-hand side of (\ref{3.1}). We want to show 
\begin{equation}
\omega = \lim_{n \rightarrow \infty} [a_1, a_2, \ldots, a_n]_m, \quad \omega \in \Omega. \label{3.22}
\end{equation}
To this end, define integer-valued functions $p_n(\omega)$ and $q_n(\omega)$, for $n \in \mathbb{N}_+$, by 
\begin{eqnarray}
p_n(\omega) &=& m^{a_n}p_{n-1}(\omega) + (m-1)m^{a_{n-1}}p_{n-2}(\omega), \quad n \geq 2, \label{3.15} \\
q_n(\omega) &=& m^{a_n}q_{n-1}(\omega) + (m-1)m^{a_{n-1}}q_{n-2}(\omega), \quad n \geq 1, \label{3.16}
\end{eqnarray}
with $p_0(\omega)=0$, $q_0(\omega) = 1$, $p_1(\omega)=1$, $q_{-1}(\omega)=0$ and $a_0 \equiv 0$. 

Now, it is easy to prove by induction that for any $n \in \mathbb{N}_+$ we have
\begin{equation}
p_n(\omega)q_{n-1}(\omega) - p_{n-1}(\omega)q_n(\omega) = (-1)^{n-1}(m-1)^{n-1}m^{a_1 + \ldots + a_{n-1}}, \label{3.17}
\end{equation}
and
\begin{equation}
\frac{m^{-a_1}}{\displaystyle 1 + \frac{(m-1)m^{-a_2}}{\displaystyle 1+ \ddots + \frac{(m-1)m^{-a_n}}{1+(m-1)t}}} = \frac{p_n(\omega) + (m-1)tm^{a_n}p_{n-1}(\omega)}{q_n(\omega) + (m-1)tm^{a_n}q_{n-1}(\omega)}, \label{3.18}
\end{equation} 
with $0 \leq t \leq 1$. 

It follows from the definitions of $\tau_m$ and $a_n$ that for any $\omega \in \Omega$ we have
\begin{equation}
\tau_m^{n-1}(\omega) = \frac{m^{-a_n}}{1+(m-1)\tau_m^n(\omega)}, \quad n \in \mathbb{N}_+, \label{3.9}
\end{equation}
hence
\begin{equation}
\omega = \frac{m^{-a_1}}{\displaystyle 1 + \frac{(m-1)m^{-a_2}}{\displaystyle 1 + \ddots+ \frac{(m-1)m^{-a_n}}{\displaystyle 1 + (m-1)\tau_m^n(\omega)}}}, \quad n \in \mathbb{N}_+. \label{3.10}
\end{equation}

By combining (\ref{3.10}), (\ref{3.15}) and (\ref{3.16}) we have
\begin{equation}
\omega = \frac{p_n(\omega) + (m-1)\tau^n_m(\omega)m^{a_n}p_{n-1}(\omega)}{q_n(\omega) + (m-1)\tau^n_m(\omega)m^{a_n}q_{n-1}(\omega)}, \quad \omega \in \Omega, \ n \in \mathbb{N}_+. \label{3.19}
\end{equation}

Taking $\tau^n_m(\omega)=0$ in (\ref{3.19}) gives
\begin{equation}
[a_1, a_2, \ldots, a_n]_m = \frac{p_n(\omega)}{q_n(\omega)}. \label{3.14}
\end{equation}

Now, using (\ref{3.17}), (\ref{3.19}) and (\ref{3.14}), for any $\omega \in \Omega$ we obtain 
\begin{equation}
\left| \omega - \frac{p_n(\omega)}{q_n(\omega)} \right| = \frac{(m-1)^n \tau^n_m(\omega)m^{a_1+\ldots+a_n}}{q_n(\omega)\left(q_n(\omega)+(m-1)\tau^n_m(\omega)m^{a_n}q_{n-1}(\omega)\right)}, n \in \mathbb{N}_+. \label{3.23}
\end{equation}
Note that this equation measure the difference between $\omega \in \Omega$ and its convergent and is the key ingredient of the following estimate.
\begin{lemma}
For any $\omega \in \Omega$ we have 
\begin{equation}
\left| \omega - \frac{p_n(\omega)}{q_n(\omega)} \right| \leq \left( \frac{m-1}{m} \right)^n, \quad n \in \mathbb{N}_+. \label{3.21}
\end{equation}
\end{lemma}
\noindent\textbf{Proof.}
By applying $\tau^{n}_m(\omega) \leq 1$ to (\ref{3.23}), we have 
\begin{equation}
\left| \omega - \frac{p_n(\omega)}{q_n(\omega)} \right| \leq \frac{(m-1)^n m^{a_1+\ldots+a_n}}{q_n(\omega)\left(q_n(\omega)+(m-1) m^{a_n}q_{n-1}(\omega)\right)}, \quad n \in \mathbb{N}_+. \label{3.23'}
\end{equation}
Let 
\begin{equation}
t_n := \frac{(m-1)^n m^{a_1+\ldots+a_n}}{q_n(\omega)\left(q_n(\omega)+(m-1) m^{a_n}q_{n-1}(\omega)\right)}, \quad n \in \mathbb{N}_+. \label{3.24'}
\end{equation}
From (\ref{3.16}), we have that $q_n(\omega)+(m-1) m^{a_n}q_{n-1}(\omega) \geq m \cdot m ^{a_n} q_{n-1}(\omega)$, i.e, $q_n(\omega) \geq m ^{a_n} q_{n-1}(\omega)$. Thus, by (\ref{3.24'}) and since $ q_n(\omega) \geq q_{n-1}(\omega)+(m-1) m^{a_{n-1}}q_{n-2}(\omega)$, we have
\begin{eqnarray}
t_n & \leq & \frac{m-1}{m} \left( \frac{(m-1)^{n-1} m^{a_1+\ldots+a_{n-1}}}{q_n(\omega)q_{n-1}(\omega)} \right) \nonumber \\ 
   & \leq & \frac{m-1}{m} \left(\frac{(m-1)^{n-1}m^{a_1+\ldots+a_{n-1}}}{q_{n-1}(\omega)\left(q_{n-1}(\omega)+(m-1) m^{a_{n-1}}q_{n-2}(\omega)\right)} \right) \nonumber \\
   & = & \frac{m-1}{m} t_{n-1}. \label{3.25'}
\end{eqnarray}
Now, by direct computation, we have 
\[
t_1 \leq \frac{m-1}{m} m^{-a_1} \leq \frac{m-1}{m} 
\]
and (\ref{3.25'}) shows that $t_n \leq \left( \frac{m-1}{m} \right)^n$, i.e., (\ref{3.21}).
\hfill $\Box$

Finally, (\ref{3.22}) follows from (\ref{3.21}), as $\frac{m-1}{m} < 1$. 

\subsection{Approximation result}

Diophantine approximation (see, e.g., \cite{Hincin-1964}) deals with the approximation of real numbers by rational numbers. Before we give the corresponding approximation result, we define the \textit{cylinder} (or \textit{fundamental interval}) of rank $n$, $I_m\left(i^{(n)}\right)$,
and show that any $I_m\left(i^{(n)}\right)$ is the set of irrationals from a certain open interval with rational endpoints.

For any $n \in \mathbb{N}_+$ and $i^{(n)}=(i_1, \ldots, i_n) \in \mathbb{N}^n$ we will say that 
\begin{equation}
I_m\left(i^{(n)}\right) = \left\{\omega \in \Omega: a_k(\omega) = i_k, 1 \leq k \leq n \right\} \label{3.26}
\end{equation}
is the \textit{fundamental interval of rank} $n$ and make the convention that $I_m\left(i^{(0)}\right) = \Omega$. 

For example, for any $i \in \mathbb{N}$ we have
\begin{equation}
I_m\left(i\right) = \left\{\omega \in \Omega: a_1(\omega) = i \right\} = \Omega \cap \left( m^{-(i+1)}, m^{-i} \right). \label{3.27}
\end{equation}
We will write $I_m(a_1, \ldots, a_n) = I_m\left(a^{(n)}\right)$, $n \in \mathbb{N}_+$. If $n \geq 2$ and $i_n \in \mathbb{N}$, then we have
\[
I_m(a_1, \ldots, a_n) = I_m\left(i^{(n)}\right).
\]
From the definition of $\tau_m$ and (\ref{3.19}) we have 
\begin{equation}
I_m\left(a^{(n)}\right) = \Omega \cap \left(u\left(a^{(n)}\right), v\left(a^{(n)}\right)\right), \label{3.28}
\end{equation}
where 
\begin{equation}
u\left(a^{(n)}\right)=\left\{
\begin{array}{lll}
	\displaystyle\frac{p_n(\omega)+(m-1)m^{a_n}p_{n-1}(\omega)}{q_n(\omega)+(m-1)m^{a_n}q_{n-1}(\omega)}, & \quad \mbox{if $n$ is odd} \\
	\displaystyle\frac{p_n(\omega)}{q_n(\omega)}, & \quad \mbox{if $n$ is even} \\ 
\end{array}
\right. \label{3.29}
\end{equation}
and
\begin{equation}
v\left(a^{(n)}\right)=\left\{
\begin{array}{lll}
	\displaystyle\frac{p_n(\omega)}{q_n(\omega)}, & \quad \mbox{if $n$ is odd} \\
	\displaystyle\frac{p_n(\omega)+(m-1)m^{a_n}p_{n-1}(\omega)}{q_n(\omega)+(m-1)m^{a_n}q_{n-1}(\omega)}, & \quad \mbox{if $n$ is even}.
\end{array}
\right. \label{3.30}
\end{equation}
Now, using (\ref{3.17}), a direct computation shows that 
\begin{equation}
\lambda\left(I\left(a^{(n)}\right)\right) = \frac{(m-1)^n m^{a_1+\ldots+a_n}}{q_n(\omega)\left(q_n(\omega)+(m-1) m^{a_n}q_{n-1}(\omega)\right)} \label{3.33}
\end{equation}
and from (\ref{3.23}) and (\ref{3.21}) we have that 
\begin{equation}
\lambda\left(I\left(a^{(n)}\right)\right) \leq \left( \frac{m-1}{m} \right)^n. \label{3.32'}
\end{equation}

We now give a Legendre-type result for these continued fraction expansions. First we define the \textit{approximation coefficient} $\Theta_m := \Theta_m (\omega)$ by
\[
\Theta_m := q^2\left|\omega-\frac{p_n}{q_n}\right|, \quad n \in \mathbb{N}_+
\]
where $\frac{p_n}{q_n}$ is the $n$th continued fraction convergent of $\omega \in \Omega$. The approximation coefficient gives a numerical indication of the quality of the approximation. 
\begin{proposition} \label{prop.2.1}
For $\omega \in \Omega$ and $p/q$ be a rational number with $p<q$, $q>0$ and $\mathrm{g.c.d.}(p, q)=1$. Let 
\[
\frac{p}{q} = [i_1, \ldots, i_n]_m, \quad \frac{p_{n-1}}{q_{n-1}}=[i_1, \ldots, i_{n-1}]_m
\]
with $p_0=0$ and $q_0=1$, where the length $n = n(p/q) \in \mathbb{N}_+$ of the continued fraction expansion of \mbox{ }$p/q$ is chosen in such a way that it is even if \mbox{ }$p/q<\omega$ and odd otherwise. Then 
\[
\Theta_m < \frac{(m-1)^n m^{i_1+\ldots+i_n}q}{q+(m-1)m^{i_n}q_{n-1}} \quad \mbox{if and only if} \quad  \frac{p}{q} \quad \mbox {is a convergent of} \quad \omega.
\]
\end{proposition}

\noindent \textbf{Proof.} 
If $p/q$ is a convergent of $\omega$, then by (\ref{3.23}) we have
\[
\Theta_m = q^2 \left|\omega - \frac{p}{q}\right| = \frac{(m-1)^n\tau^n_m(\omega)m^{i_1+\ldots+i_n}q}{q+(m-1)\tau^n_m(\omega)m^{i_n}q_{n-1}(\omega)} \leq \frac{(m-1)^nm^{i_1+\ldots+i_n}q}{q+(m-1)m^{i_n}q_{n-1}}.
\]
Conversely, if $\Theta_m < \displaystyle\frac{(m-1)^nm^{i_1+\ldots+i_n}q}{q+(m-1)m^{i_n}q_{n-1}}$, then
\[
q\left|\omega-\frac{p}{q}\right| < \frac{(m-1)^nm^{i_1+\ldots+i_n}}{q+(m-1)m^{i_n}q_{n-1}}.
\]
Assuming that $n$ is even, then $\omega > \displaystyle\frac{p}{q}$ and we have $\omega-\displaystyle\frac{p}{q} < \displaystyle\frac{(m-1)^nm^{i_1+\ldots+i_n}}{q(q+(m-1)m^{i_n}q_{n-1})}$. Thus,
\[
\frac{p}{q}< \omega < \frac{p}{q} + \frac{(m-1)^nm^{i_1+\ldots+i_n}}{q(q+(m-1)m^{i_n}q_{n-1})} = \frac{p+(m-1)m^{i_n}p_{n-1}}{q+(m-1)m^{i_n}q_{n-1}}.
\]
Hence, $\omega \in I_m\left(i^{(n)}\right)$, i.e., $\displaystyle\frac{p}{q}=[i_1, \ldots, i_n]_m$ is a convergent of $\omega$. The case when $n$ is an odd is treated similarly. 
\hfill $\Box$

\subsection{The probability structure of $(a_n)_{n \in \mathbb{N}_+}$ under the $\lambda$}

We start by deriving the so-called Brod\'en-Borel-L\'evy formula (see, e.g., \cite{r10, r11}) for these type of expansions. First, define $s_n$, $n \in \mathbb{N}_+$, by 
\begin{equation}
s_n = m^{-a_n}\frac{q_n}{q_{n-1}}-1, \quad s_1=0, \label{3.34}
\end{equation}
where $m \geq 2$ and $a_n$, $q_n$ are defined in (\ref{1.3}) and (\ref{3.16}), respectively.

Next, (\ref{3.16}) implies that
\begin{equation}
s_n = \frac{(m-1)m^{-a_n}}{1+s_{n-1}}, \quad n \geq 2, \label{3.35}
\end{equation}
hence 
\begin{equation}
s_n = \frac{(m-1)m^{-a_n}}{\displaystyle 1+\frac{(m-1)m^{-a_{n-1}}}{\displaystyle 1 +\ddots+ \frac{(m-1)m^{-a_3}}{\displaystyle 1 + (m-1)m^{-a_2}}}} = (m-1)[a_n, a_{n-1}, \ldots, a_2, \infty]_m, \label{3.36}
\end{equation} 
for $n \geq 2$.
\begin{proposition} [Brod\'en-Borel-L\'evy formula type] 
For any $n \in \mathbb{N}_+$ we have
\begin{equation}
\lambda\left(\tau^n_m < x |a_1,\ldots, a_n \right) = \frac{(s_n + m)x}{(s_n+(m-1)x+1)}, \quad x \in I, \label{3.37}
\end{equation}
where $s_n$ is defined by (\ref{3.34}) or (\ref{3.35}).
\end{proposition}
\noindent\textbf{Proof.} As we know, for any $n \in \mathbb{N}_+$ and $x \in I$, we have
\[
\lambda\left(\tau^n_m < x |a_1,\ldots, a_n \right) = \frac{\lambda\left(\left(\tau^n_m < x\right) \cap I_m(a_1,\ldots, a_n) \right)}{\lambda\left(I_m(a_1,\ldots, a_n)\right)}. 
\]
From (\ref{3.19}) and (\ref{3.28}) we have
\begin{eqnarray}
\lambda\left(\left(\tau^n_m < x \right) \cap I(a_1,\ldots, a_n)\right) &=& \left|\frac{p_n}{q_n} - \frac{p_n+(m-1)xm^{a_n}p_{n-1}}{q_n+(m-1)xm^{a_n}q_{n-1}}\right| \nonumber \\
&=& \frac{(m-1)^nxm^{a_1+\ldots+ a_n}}{q_n\left(q_n+(m-1)xm^{a_n}q_{n-1}\right)}. \nonumber
\end{eqnarray}
Hence, from (\ref{3.33}) we have 
\begin{eqnarray*}
\lambda\left(\tau^n_m < x |a_1, \ldots, a_n \right) &=& \frac{\lambda\left(\left(\tau^n_m < x\right) \cap I_m(a_1,\ldots, a_n) \right)}{\lambda\left(I_m(a_1,\ldots, a_n)\right)} \\  
&=& \frac{x\left(q_n+(m-1)m^{a_n}q_{n-1}\right)}{q_n\left(q_n+(m-1)xm^{a_n}q_{n-1}\right)} \\  
&=& \frac{(s_n+m)x}{s_n+(m-1)x+1}, 
\end{eqnarray*}
for any $n \in \mathbb{N}_+$ and $x \in I$. 
\hfill $\Box$

The Brod\'en-Borel-L\'evy formula allows us to determine the probability structure of $(a_n)_{n \in \mathbb{N}_+}$ under $\lambda$. 
\begin{proposition}
For any $i \in \mathbb{N}$ and $n \in \mathbb{N}_+$ we have
\begin{equation}
\lambda(a_1=i) = (m-1)m^{-(i+1)} \label{3.38}
\end{equation}
and
\begin{equation}
\lambda\left(a_{n+1}=i |a_1,\ldots, a_n \right) = P^i_m(s_n), \label{3.39}
\end{equation}
where 
\begin{equation}
P^i_m(x) = \frac{(m-1)m^{-(i+1)}(x+1)(x+m)}{(x+(m-1)m^{-i}+1)(x+(m-1)m^{-(i+1)}+1)}. \label{3.40}
\end{equation}
\end{proposition}
\noindent\textbf{Proof.}
As shown above, we have
\[
\left\{\omega \in \Omega: a_1(\omega) = i\right\} = \Omega \cap \left(m^{-(i+1)}, m^{-i}\right). 
\]
Thus,
\[
\lambda(a_1=i) = \left|m^{-(i+1)}- m^{-i}\right| = (m-1)m^{-(i+1)}.
\]
From (\ref{3.9}), we have that 
\[
\tau_m^n(\omega) = [a_{n+1}, a_{n+2}, \ldots]_m, \quad n \in \mathbb{N}_+, \quad \omega \in \Omega
\]
and so we have
\begin{eqnarray*}
\lambda\left(a_{n+1}=i |a_1,\ldots, a_n \right) & = & \lambda\left( \tau^n_m \in \left(m^{-(i+1)}, m^{-i}\right] | a_1,\ldots, a_n \right) \\
&=& \frac{(s_n+m)m^{-i}}{s_n+(m-1)m^{-i}+1} - \frac{(s_n+m)m^{-(i+1)}}{s_n+(m-1)m^{-(i+1)}+1} \\
&=& \frac{(m-1)m^{-(i+1)}(s_n+1)(s_n+m)}{(s_n+(m-1)m^{-i}+1)(s_n+(m-1)m^{-(i+1)}+1)} \\
&=& P^i_m(s_n).  
\end{eqnarray*}
\hfill $\Box$

Hence, the sequence $(s_n)_{n \in \mathbb{N}_+}$ with $s_1=0$ is a homogeneous $I$-valued Markov chain on $\left(I, {\mathcal{B}}_I, \lambda \right)$ with the following transition mechanism: from state $s \in I \setminus \Omega$, $s \geq 1$ the only possible one-step transitions are those to states $m^{-i}/(1+(m-1)s)$, $i \in \mathbb{N}$, with corresponding probabilities $P^i_m(s)$, $i \in \mathbb{N}$.  

\subsection{The invariant measure of $\tau_m$}

In this subsection we will give the explicit form of the invariant probability measure $\gamma_m$ of the transformation $\tau_m$, i.e., $\gamma_m (A) = \gamma_m \left(\tau^{-1}_m(A)\right)$, $A \in {\mathcal{B}}_I$. 

Let ${\mathcal{B}}_I$ denote the $\sigma$-algebra of Borel subsets of $I$. The metric point of view in studying the sequence $(a_n)_{n \in \mathbb{N}_+}$ is to consider that the $a_n$, $n \in \mathbb{N}_+$, are non-negative integer-valued random variables which are defined almost surely on $\left(I, {\mathcal{B}}_I\right)$ with respect to any probability measure on ${\mathcal{B}}_I$ that assign probability $0$ to the set $I \setminus \Omega$ of rationals in $I$. Such a measure is Lebesgue measure $\lambda$.

Another measure on ${\mathcal{B}}_I$ more important than Lebesgue measure, that assign probability $0$ to the set of rationals in $I$, is the \textit{invariant probability measure} $\gamma_m$ of the transformation $\tau_m$.

\begin{proposition} \label{prop.3.3.1}
The invariant probability density $\rho_m$ of the transformation $\tau_m$ is given by
\begin{equation}
\rho_m(x)=\frac{1}{((m-1)x+1)((m-1)x+m)}, \quad x \in I, \label{3.41}
\end{equation}
with the normalizing factor $k_m = \frac{(m-1)^2}{\log \left(m^2/(2m-1)\right)}$.
\end{proposition}
\noindent \textbf{Proof.} See Appendix.

\noindent Hence 
\begin{equation}
\gamma_m (A) = k_m \int_{A} \frac{dx}{((m-1)x+1)((m-1)x+m)}, \quad A \in {\mathcal{B}}_I. \label{3.44}
\end{equation}

The normalization constant $k_m$ defined above is chosen so that $\gamma_m ([0, 1]) = 1$. 

\section{The natural extension of $\tau_m$ and extended random variables}

By its very definition, the sequence $(a_n)_{n \in \mathbb{N}_+}$ in (\ref{1.3}) and (\ref{1.4}) is strictly stationary under $\gamma_m$. As such, there should exist a doubly infinite version of it, say $\overline{a}_l$, $l \in \mathbb{Z}:=\left\{\ldots, -1, 0, 1, \ldots\right\}$, defined on a richer probability space. It appears that this doubly infinite version can be effectively constructed on $(I^2, {\mathcal{B}}^2_I, \overline{\gamma}_m)$, where $\overline{\gamma}_m$ is the so-called extended measure which expresion is given below.

\subsection{Definition and basic properties}

For $\tau_m$ in (\ref{3.5}), the \textit{natural extension} $\overline{\tau}_m$ of $\tau_m$ \cite{r20} is the transformation of $\left[0, 1\right) \times I$ defined by
\begin{equation}
\overline{\tau}_m(x,y) = \left( \tau_m(x), \frac{m^{-a_1(x)}}{(m-1)y+1} \right), \mbox{ } (x, y) \in \left[0, 1\right) \times I. \label{3.45}
\end{equation}
This is a one-to-one transformation of $\Omega^2$ with the inverse
\begin{equation}
\overline{\tau}^{-1}_m(\omega, \theta) = \left(\frac{m^{-a_1(\theta)}}{(m-1)\omega+1}, \tau_m(\theta) \right), \mbox{ } (\omega, \theta) \in \Omega^2. \label{3.46}
\end{equation}
It is easy to check that for $n \geq 2$ we have 
\begin{equation}
\overline{\tau}^n_m(\omega, \theta) = \left(\tau^n_m(\omega), \left[a_n(\omega), a_{n-1}(\omega), \ldots, a_2(\omega), a_1(\omega) + \frac{\log(1+(m-1)\theta)}{\log m} \right]_m \right), \label{3.47}
\end{equation}
and
\begin{equation}
\overline{\tau}^{-n}_m(\omega, \theta)=\left(\left[a_n(\theta), a_{n-1}(\theta), \ldots, a_2(\theta), a_1(\theta)+\frac{\log(1+(m-1)\omega)}{\log m} \right]_m, \tau^{n}_m(\theta) \right). \label{3.48}
\end{equation}
Now, define the \textit{extended measure} $\overline{\gamma}_m$ on ${\mathcal{B}}^2_I$ as
\begin{equation}
\overline{\gamma}_m (B) = k_m \int\!\!\!\int_{B} \frac{dxdy}{((m-1)(x+y)+1)^2} , \quad B \in {\mathcal{B}}^2_I. \label{3.49}
\end{equation}
A simple calculus show us that 
\begin{equation}
\overline{\gamma}_m(A \times I) = \overline{\gamma}_m(I \times A) = \gamma_m(A), \mbox{ } A \in {\mathcal{B}}_I. \label{3.50}
\end{equation}
The result below shows that $\overline{\gamma}_m$ plays with respect to $\overline{\tau}_m$ the part played by $\gamma_m$ with respect to $\tau_m$.

\begin{proposition} \label{prop.3.1}
The extended measure $\overline{\gamma}_m$ is preserved by $\overline{\tau}_m$. 
\end{proposition}
\noindent\textbf{Proof.} See Appendix.

\subsection{Extended random variables}

Define extended incomplete quotients $\overline{a}_l$, $l \in \mathbb{Z}$, on $\Omega^2$ by
\[
\overline{a}_{l+1}(\omega, \theta) = \overline{a}_1 \left(\overline{\tau}_m^l (\omega, \theta) \right), \quad l \in \mathbb{Z},
\]
with 
\[
\overline{a}_1(\omega, \theta) = a_1 (\omega), \quad (\omega, \theta) \in \Omega^2.
\]
By (\ref{3.47}) and (\ref{3.48}) we have
\[
\overline{a}_n(\omega, \theta) = a_n(\omega), \mbox{ } \overline{a}_0(\omega, \theta) = a_1(\theta), \mbox{ } \overline{a}_{-n}(\omega, \theta) = a_{n+1}(\theta), \quad n \in \mathbb{N}_+, \ (\omega, \theta) \in \Omega^2. 
\]

\begin{remark} 
Since $\overline{\tau}_m$ preserves $\overline{\gamma}_m$, the doubly infinite sequence $\left(\overline{a}_l\right)_{l \in \mathbb{Z}}$, is strictly stationary under $\overline{\gamma}_m$.
\end{remark}

\begin{theorem} \label{th.3.3.3}
For any $x \in I$ we have
\begin{equation}
\overline{\gamma}_m ( [0, x] \times I \left. \right| \overline{a}_0, \overline{a}_{-1}, \ldots ) = \frac{((m-1)a + m)x}{(m-1)(x+a) + 1} \quad \overline{\gamma}_m \mbox{-}\mathrm{a.s.}, \label{3.55}
\end{equation}
where $a = [\overline{a}_0, \overline{a}_{-1}, \ldots]_m$. 
\end{theorem}
\noindent \textbf{Proof.} Let $I_{m,n}$ denote the fundamental interval $I_m(\overline{a}_0, \overline{a}_{-1}, \ldots, \overline{a}_{-n})$, $n \in \mathbb{N}$. We have
\[
\overline{\gamma}_m ( [0, x] \times I \left. \right| \overline{a}_0, \overline{a}_{-1}, \ldots ) = \lim_{n \rightarrow \infty} \overline{\gamma}_m ( [0, x] \times I \left. \right| \overline{a}_0, \ldots, \overline{a}_{-n} ) \quad \overline{\gamma}_m \mbox{-a.s.}
\]
and
\begin{eqnarray*}
\overline{\gamma}_m ( [0, x] \times I \left. \right| \overline{a}_0, \ldots, \overline{a}_{-n} ) = \frac{\overline{\gamma}_m ([0, x] \times I_{m,n})}{\overline{\gamma}_m (I \times I_{m,n})} \qquad \qquad \qquad \qquad \qquad \\
 =  \frac{k_m \displaystyle\int_{I_{m,n}}dy\displaystyle\int^x_0{\frac{du}{((m-1)(u+y)+1)^2}}}{\gamma_m(I_{m,n})} \qquad \qquad \qquad \quad \ \\
 =  \frac{1}{\gamma_m(I_{m,n})} k_m \int_{I_{m,n}} \frac{x}{((m-1)(x+y)+1)((m-1)y+1)}dy \\
 =  \frac{1}{\gamma_m(I_{m,n})} \int_{I_{m,n}} \frac{x((m-1)y+m)}{(m-1)(x+y)+1} \gamma_m(dy) \qquad \qquad \quad \ \ \\
 =  \frac{x((m-1)y_n+m)}{(m-1)(x+y_n)+1}, \qquad \qquad \qquad \qquad \qquad \qquad \quad 
\end{eqnarray*}
for some $y_n \in I_{m,n}$. Since 
\begin{equation}
\lim_{n \rightarrow \infty} y_n =[\overline{a}_0, \overline{a}_{-1}, \ldots]_m = a, \label{3.56}
\end{equation}
the proof is complete.
\hfill $\Box$

The stochastic property of $(\overline{a}_l)_{l \in \mathbb{Z}}$ under $\overline{\gamma}_m$ is given by the following corollary of Theorem \ref{th.3.3.3}. 
\begin{corollary}
For any $i \in \mathbb{N}$ we have
\[
\overline{\gamma}_m (\left.\overline{a}_1 = i\right| \overline{a}_0, \overline{a}_{-1}, \ldots) = P^i_m((m-1)a) \quad \overline{\gamma}_m \mbox{-}\mathrm{a.s.},
\]
where $a = [\overline{a}_0, \overline{a}_{-1}, \ldots]_m$.
\end{corollary}
\noindent \textbf{Proof.} Let us denote by $I_{m,n}$ the fundamental interval $I_m(\overline{a}_0, \overline{a}_{-1}, \ldots, \overline{a}_{-n})$, $n \in \mathbb{N}$. We have
\[
(\overline{a}_1 = i) = \left( m^{-(i+1)}, m^{-i} \right] \times [0, 1)
\]
and
\[
\overline{\gamma}_m (\left.\overline{a}_1 = i\right| \overline{a}_0, \overline{a}_{-1}, \ldots) = \lim_{n \rightarrow \infty} \overline{\gamma}_m (\left.\overline{a}_1 = i\right| I_{m,n}).
\]
Now 
\begin{eqnarray}
\overline{\gamma}_m \left( \left.\left( m^{-(i+1)}, m^{-i} \right) \times [0, 1)\right| I_{m,n}\right) &=& \frac{\overline{\gamma}_m \left(\left( m^{-(i+1)}, m^{-i} \right) \times I_{m,n}\right)}{\overline{\gamma}_m (I \times I_{m,n})} \nonumber \\
& = & \frac{1}{\gamma_m (I_{m,n})} \int_{I_n} P^i_m((m-1)y) \gamma_m(dy) \nonumber \\
& = & P^i_m((m-1)y_n), \nonumber
\end{eqnarray}
for some $y_n \in I_{m,n}$. From (\ref{3.56}) the proof is complete. 
\hfill $\Box$

\begin{remark}
The strict stationarity of $\left(\overline{a}_l\right)_{l \in \mathbb{Z}}$, under $\overline{\gamma}_m$ implies that 
\[
\overline{\gamma}_m (\left.\overline{a}_{l+1} = i\right| \overline{a}_l, \overline{a}_{l-1}, \ldots) = P^i_m((m-1)a) \quad \overline{\gamma}_m-a.s. 
\]
for any $i \in \mathbb{N}$ and $l \in \mathbb{Z}$, where $a = [\overline{a}_l, \overline{a}_{l-1}, \ldots]_m$. The last equation emphasizes that $\left(\overline{a}_l\right)_{l \in \mathbb{Z}}$ is a chain of infinite order in the theory of dependence with complete connections (see \cite{r10}, Section 5.5).
\end{remark}

Motivated by Theorem \ref{th.3.3.3} we shall consider the family of (conditional) probability measures $\left(\gamma^a_m\right)_a$ on ${\mathcal{B}}_I$ defined by their distribution functions
\begin{equation}
\gamma^a_m ([0, x]) = \frac{((m-1)a +m)x}{(m-1)(x+a) + 1}, \ x \in I, \ a \geq 0. \label{3.57}
\end{equation}
Note that the limit case $a = \infty$ is $\gamma^\infty_m = \lambda$. 

For any $a \geq 0$ put $s^a_0 = a$ and
\begin{equation}
s^a_n = \frac{(m-1)m^{-a_n}}{1+s^a_{n-1}}, \ n \in \mathbb{N}_+. \label{3.58}
\end{equation}
For $a \geq 0$ we have 
\[
s^a_1 = \frac{(m-1)m^{-a_1}}{1+a}
\]
and
\[
s^a_n = (m-1)\left[a_n, \ldots, a_2, a_1+\frac{\log(a+1)}{\log m}\right]_m, \quad n \geq 2. 
\]
Then $\left(s^a_n\right)_{n \in \textbf{N}_+}$ is a $I \cup \left\{a\right\} $ - valued Markov chain on $(I, {\mathcal{B}}_I, \gamma^a_m)$ which starts from $s^a_0 = a\geq0$ and has the following transition mechanism: from state $s \in I \cup \left\{a\right\}$ the possible transitions are to any state $m^{-i}/((m-1)s+1)$ with the corresponding transition probability $P^i_m((m-1)s)$, $i \in \mathbb{N}$. 

Now, it is easy to check by induction that
\begin{equation}
s^a_n = m^{-a_n} \frac{(m-1)p_n+(a+1)q_n}{(m-1)p_{n-1}+(a+1)q_{n-1}} - 1, \label{3.72'}
\end{equation}
for any $n \in \mathbb{N}_+$ and $a \geq 0$.

Thus, a simple calculation shows that for any $n \in \mathbb{N}_+$ we have
\begin{eqnarray*}
\gamma^a_m \left( \left.\tau^n_m < x\right| a_1, \ldots, a_n \right) = \frac{\gamma^a_m \left(\left(\tau^n_m < x\right) \cap I_m \left(a^{(n)}\right)\right)}{\gamma^a_m \left(I_m(a^{(n)})\right)} \qquad \qquad \qquad \qquad \qquad \\
= \frac{x((m-1)((m-1)p_n+(a+1)q_n)+m^{a_n}((m-1)p_{n-1}+(a+1)q_{n-1}))}{(m-1)((m-1)p_n+(a+1)q_n)+xm^{a_n}((m-1)p_{n-1}+(a+1)q_{n-1})}.
\end{eqnarray*}
By (\ref{3.72'}) for any $ n \in \mathbb{N}_+$ we have 
\begin{equation}
\gamma^a_m \left( \left.\tau^n_m < x\right| a_1, \ldots, a_n \right) = \frac{((m-1)s^a_n +m)x}{(m-1)(x+s^a_n) + 1}, \mbox{ } a \geq 0, x \in I. \label{3.59}
\end{equation}
The last equation is the generalization of the Brod\'en-Borel-L\'evy formula from section 2.3.

\section{The Perron-Frobenius operator of $\tau_m$ under $\gamma_m$}

In this section we derive and study the associated Perron-Frobenius operator of $\tau_m$ under the invariant measure $\gamma_m$.

Let $\mu$ be a probability measure on $\left(I, {\mathcal{B}}_I\right)$ such that $\mu\left(\tau^{-1}_m(A)\right) = 0$ whenever $\mu(A) = 0$, $A \in {\mathcal{B}}_I$, where the transformation $\tau_m$ is defined in (\ref{3.5}). In particular, this condition is satisfied if $\tau_m$ is $\mu$-preserving, that is, $\mu \tau^{-1}_m = \mu$. It is known from previous section, that the Perron-Frobenius operator $P_{\mu}$ of $\tau_m$ under $\mu$ is defined as the bounded linear operator on $L^1_{\mu}=\left\{f: I \rightarrow \mathbb{C} | \int_I \left|f\right|d\mu < \infty \right\}$ which takes $f \in L^1_{\mu}$ into $P_{\mu}f \in L^1_{\mu}$ with
\begin{equation}
\int_{A}P_{\mu}f d\mu = \int_{\tau^{-1}_m(A)}f d\mu, \quad A \in {\mathcal{B}}_I. \label{3.60}
\end{equation}
In particular, the Perron-Frobenius operator $P_{\lambda}$ of $\tau_m$ under the Lebesgue measure $\lambda$ is (see \cite{r1}, p.86) 
\begin{equation}
P_{\lambda}f(x)=\frac{d}{dx}\int_{\tau^{-1}_m([0,x])}f d\lambda=\sum_{t \in \tau^{-1}_m(x)}\frac{f(t)}{\left|\tau'_m(t)\right|} \mbox{  } \mbox{  } \mbox{a.e. in } I. \label{3.61}
\end{equation}
The following results will be proved in the Appendix.

The following Proposition gives the expression of the Perron-Frobenius operator of $\tau_m$ under the invariant measure $\gamma_m$ (\ref{3.62}) and under a probability measure which is absolutely continuous with respect to the Lebesgue measure (\ref{3.69}). Also, we derive the asymptotic behaviour of this operator (\ref{3.75}). 

\begin{proposition} \label{prop.4.1}
\begin{enumerate}
	\item[(i)] The Perron-Frobenius operator $U_m := P_{\gamma_m}$ of $\tau_m$ under $\gamma_m$ is given a.e. in $I$ by the equation 
\begin{equation}
U_mf(x) = \sum_{i \in \mathbb{N}}P^i_m((m-1)x)f(u^i_m(x)), \mbox{ } f \in L^1_{\gamma_m}, \label{3.62}
\end{equation}
where $P^i_m$ is defined in $(\ref{3.40})$ and $u^i_m(x)$ is given by the equation
\begin{equation}
u^i_m(x) = \frac{m^{-i}}{(m-1)x+1}, \mbox{ } x \in I. \label{3.63}
\end{equation}

\item[(ii)] Let $\mu$ be a probability measure on ${\mathcal{B}}_I$. Assume that $\mu$ is absolutely continuous with respect to $\lambda$ (and denote $\mu \ll \lambda$, i.e., if $\mu(A) = 0$ for every set $A$ with $\lambda(A) = 0$) and let $h = d\mu / d \lambda$  \mbox{a.e. in } I. Then:
\begin{enumerate}
	\item[(a)] the Perron-Frobenius operator $P_{\mu}$ of $\tau_m$ under $\mu$ is given a.e. in $I$ by the equation 
	\begin{eqnarray}
P_{\mu} f(x) &=& \frac{1}{h(x)} \sum_{i \in \mathbb{N}} \frac{h(u^i_m(x))}{((m-1)x+1)^2} (m-1)m^{-i}f(u^i_m(x)) \label{3.68} \\
             &=& \frac{U_m g(x)}{((m-1)x+1)((m-1)x+m)h(x)}, \ f \in L^1_{\mu}, \label{3.69}
\end{eqnarray}
where $g(x)=((m-1)x+1)((m-1)x+m)f(x)h(x)$, $x \in I$.

The powers of $P_{\mu}$ are given a.e. in $I$ and for any $f \in L^1_{\mu}$ and any $n \in \mathbb{N}_+$ by the equation 
\begin{equation}
P^n_{\mu} f(x) = \frac{U^n_m g(x)}{((m-1)x+1)((m-1)x+m)h(x)}. \label{3.70}
\end{equation}
	\item[(b)] we have 
	\begin{equation}
\mu \left(\tau^{-n}_m(A)\right) = \int_{A} \frac{U^n_mf(x)}{((m-1)x+1)((m-1)x+m)}dx, \label{3.75}
  \end{equation}
for any $n \in \mathbb{N}$ and $A \in {\mathcal{B}}_I$, where $f(x)=((m-1)x+1)((m-1)x+m)h(x)$, $x \in I$.

\end{enumerate}

\end{enumerate}
\end{proposition}

In the next Proposition the domain of $U_m$ will be successively restricted to the following Banach spaces: $BV(I)$-the linear space of all complex-valued functions of bounded variation and $B(I)$ is the collection of all bounded measurable functions $f:I \rightarrow \mathbb{C}$. The variation $\mathrm{var}_A f$ over $A \subset I$ of a function $f:I \rightarrow \mathbb{C}$ is defined as
\[
\sup \sum^{k-1}_{i=1} |f(t_i) - f(t_{i-1})|,
\]
the supremum being taken over $t_1 < \ldots < t_k$, $t_i \in A$, $1 \leq i \leq k$, and $k \geq 2$. We write simply $\mathrm{var} f$ for $\mathrm{var}_I f$.

\begin{proposition} \label{prop.4.2} 
\begin{enumerate}
\item[(i)] If $f \in BV(I)$ is a real-valued function, then
\begin{equation}
\mathrm{var}\mbox{ } U_mf \leq K_m \mathrm{var } f, \label{3.77}
\end{equation}
where $K_m = \displaystyle\frac{\left(m-1\right)\left(3m^2-3m+1\right)}{(2m-1)\left(m^2+m-1\right)}$. The constant cannot be lowered. 
\item[(ii)] The operator $U_m : B(I) \rightarrow B(I)$ is the transition operator of the Markov chain $(s^a_n)_{n \in \mathbb{N}_+}$ on $(I, {\mathcal{B}}_I, \gamma^a_m)$, for any $a \in I$, where $(s^a_n)_{n \in \mathbb{N}_+}$ and $\gamma^a_m$ are give in (\ref{3.57}) and (\ref{3.58}), respectively.
\end{enumerate}
\end{proposition} 

\section{Proof of the Gauss-Kuzmin-type theorem}

In this section we prove our main theorem. The main tool of this section is the random system with complete connections. We will first give a brief introduction to the theory of random systems with complete connections and list some of the main applications and some important properties. The general concepts presented here will be customized in the second subsection for the continued fraction expansion presented in this paper. All these concepts will be applied in subsection 5.3 to solve our main theorem.

\subsection{Random systems with complete connections}

The purpose of this subsection is to recall the definition of random systems with complete connections, and take this opportunity to inform nonspecialists a little about some applications of the theory of random systems with complete connections.

The first explicit formal definition of the concept of dependence with complete connections was given by Onicescu and Mihoc in the 1930's when studying so-called \textit{urn schemes} (see, e.g., \cite{r21}, or \cite{r13} or the Introduction in \cite{r10}). The concept of random system with complete connections was defined by Iosifescu \cite{MI-1963}.
There are many other areas where the theory of RSCC can be applied. Let us just mention a few: \textit{mathematical modelling of learning processes} (see, e.g., \cite{N-1972, r13, Karlin-1953}), \textit{chains of infinite order} (see, e.g., \cite{r7, r8}), \textit{partially observed random chains} (see, e.g., \cite{r13}), \textit{image coding} (see \cite{BE-1988}), and \textit{continued fraction expansion} (see \cite{r10}). Nowadays RSCC are called \textit{iterated functions systems with place-dependent probabilities} or simply \textit{iterated functions systems} (IFS). This terminology was introduced by Barnsley et al. in the middle of the 1980's in \cite{Barnsley-1988}. It only became fashionable in the framework of fractals and chaos but, before that, it appeared as the simplest case of a random system with complete connections and, in particular, as the Bush-Mosteller model for learning with experimenter-controlled-events [see, e.g., \cite{BE-1988, HIR-2003}]. An application of IFS to continued fractions can be found in the paper \cite{MU-2003}. 

\subsubsection{Definitions and explanations}

First, let $(W, {\mathcal W})$ and $(X, {\mathcal X})$ be two measurable spaces. A real valued function $P$ defined on $W \times {\mathcal X}$ is called a \textit{transition probability function} from $(W, {\mathcal W})$ to $(X, {\mathcal X})$ if $P(w, \cdot)$ is a probability on ${\mathcal X}$ for any $w \in W$ and $P(\cdot, A)$ is a ${\mathcal W}$-measurable function for any $A \in {\mathcal X}$.

A quadruple 
\begin{equation}
\left\{(W, {\mathcal W}), (X, {\mathcal X}), u, P\right\} \label{5.1}
\end{equation}
is named a \textit{random system with complete connections} (RSCC) if
\begin{enumerate}
	\item[(i)] $(W, {\mathcal W})$ and $(X, {\mathcal X})$ are measurable spaces;
	\item[(ii)] $u: W \times X \rightarrow W$ is a $({\mathcal W} \otimes {\mathcal X}, {\mathcal W})$-measurable function;
	\item[(iii)] $P$ is a transition probability function from $(W, {\mathcal W})$ to $(X, {\mathcal X})$.
\end{enumerate}

The definition of a RSCC can be extended to the non-homogeneous case in the sense that all the entities constituting it are allowed to depend on $t \in T$, where $T$ is either the set $\mathbb{N}$ of natural numbers or the set $\mathbb{Z}$ of integers.

The set $W$ is usually called the \textit{state space}, the set $X$ is often called the \textit{event space} and the function $u$ is often called the \textit{response-function}. We also call $u(\cdot, x) : W \rightarrow W$ a response-function. 

The interpretation of this structure is as follows. If $X$ denotes the set of possible observations and $W$ the range of possible states of the system, then $P$ induces for every state $w \in W$ the distribution $P(w, \cdot)$ of the random observation following $w$. The function $u$ represents the transition function of the system, which transforms a given state $w$ and an actual observation $x$ into a new state $u(w, x)$.

To every RSCC $\left\{(W, {\mathcal W}), (X, {\mathcal X}), u, P\right\}$  and every $w \in W$ (an arbitrary fixed element of $W$) one can generate two stochastic sequences $\left\{\xi_n\right\}_{n \in \mathbb{N}}$ and $\left\{\zeta_n\right\}_{n \in \mathbb{N}_+}$ as follows: we set $\xi_0 = w$, pick an element $\zeta_1 \in X$ using $P(\xi_0, \cdot)$, define $\xi_1 = u(\xi_0, \zeta_1)$, pick $\zeta_2$ in $X$ using $P(\xi_1, \cdot)$, define $\xi_2 = u(\xi_1, \zeta_2)$, and generally we pick $\zeta_n$ in $X$ using $P(\xi_{n-1}, \cdot)$, and define $\xi_n = u(\xi_{n-1}, \zeta_n)$. Thus, the two stochastic sequences can be described as follows:
\begin{eqnarray*}
\xi_0 = w, \ \xi_{n+1} = u(\xi_{n}, \zeta_{n+1}), \ n \geq 1, \qquad \qquad \qquad \quad \\
P(\zeta_1 \in A) = P(w, A), \ A \in {\mathcal X} \qquad \qquad \qquad \qquad \qquad  \\
P(\zeta_{n+1} \in A | \xi_n, \zeta_n, \ldots, \xi_1, \zeta_1, \xi_0) = P(\xi_n, A), \ A \in {\mathcal X}.
\end{eqnarray*}
We call the sequence $\left\{\xi_n\right\}_{n \in \mathbb{N}}$ of $W$-valued random variables the \textit{state sequence} and the sequence $\left\{\zeta_n\right\}_{n \in \mathbb{N}_+}$ of $X$-valued random variables the \textit{event sequence}. When we want to emphasize the initial point $w$, we write
\[
\xi_n = \xi_n(w) \quad \mbox{and} \quad \zeta_n = \zeta_n(w).
\]
The central issue in the theory of dependence with complete connections is the sequence $\{\zeta_n\}_{n \in \mathbb{N}_+}$ which is a stochastic process that is no longer Markovian, but a chain with complete connections (processes whose transition probabilities depend on the whole past history). 

From the definition of $\xi_n$ it is clear that the state sequence $\left\{\xi_n\right\}_{n \in \mathbb{N}}$ is a Markov chain (the so-called associated Markov chain) with transition probability function $Q$, where 
\begin{equation}
Q(w, A) = P(w, \{x \in X | u(w, x) \in A\}) \label{5.2} 
\end{equation}
with $A \in {\mathcal W}$. 

The transition operator $U : B(W, {\mathcal W}) \rightarrow B(W, {\mathcal W})$ is defined by
\begin{equation}
Uf(w) = \sum_{x \in X} P(w, x)f(u(w,x)), \quad f \in B(W, {\mathcal W}),  \label{5.3}
\end{equation}
where $B(W, {\mathcal W})$ is the Banach space of all bounded ${\mathcal W}$-measurable complex-valued functions defined on $W$.

\subsubsection{Examples of RSCCs}
In this section we shall give two examples of RSCCs which occur either in various chapters of probability theory or as a result of modelling phenomena in various fields.
\begin{example} 
The concept of a random system with complete connections may be regarded as a generalization and formalization of the notion of a stochastic learning model. Learning may be defined as an adaptive modification of behaviour in the course of repeated trials. By mathematical learning theory  we mean the body of research methods and results concerned with the conceptual representation of learning phenomena, the mathematical formulation of hypotheses about learning, and the derivation of testable theorems. The purpose of mathematical learning theory is to provide simple, quantitative descriptions of processes which are basic to behavioural modifications. 

All stochastic models for learning studied so far fit the following general theoretical scheme. The behaviour of the subject on trial $n$ is determined by its state $S_n$ (an indicator of the subject's tendencies) at the beginning of the trial. Here $S_n$ is a random variable taking values in a measurable space $(S, {\mathcal S})$. On trial $n$ an event $E_{n+1}$ occurs that results in a change of the state. Here $E_{n+1}$ is a random variable taking values in the measurable space $(E, {\mathcal E})$ and specifies those occurrences on trial $n$ that affect the subsequent behaviour. To represent the fact that the occurrence of an event affects a change of state it is necessary to consider a measurable map $v$ from $S \times E$ into $S$ and postulate that $S_{n+1} = v(S_{n}, E_{n+1})$, $n \in \mathbb{N}$. Finally assume that the probability distribution of $E_{n+1}$ given $S_n, E_n, \ldots, S_1, E_1, S_0$ depends only on the state $S_n$ and denote it by $R(S_n, \cdot)$. By a general learning model we mean the collection $\left\{(S, {\mathcal S}), (E, {\mathcal E}), v, R\right\}$ which is trivially an RSCC. Notice that in fact we only changed the notation. Various special learning models are obtained by simply particularizing $S$, $E$, $v$ and $R$ (see, e.g., \cite{r13, N-1972}). 
\hfill $\Box$
\end{example} 

\begin{example}
As we mentioned in subsection 1.1, any irrational number $y$ in the unit interval $[0,1]$ has an infinite continued fraction expansion of the form 
\[
y = \displaystyle \frac{1}{a_1(y)+\displaystyle \frac{1}{a_2(y)+\displaystyle \frac{1}{a_3(y) + \ddots}}},
\] 
where the $a_n(y)$, $n \in \mathbb{N}_+$, are natural numbers. Define $(s_n)_{n \in \mathbb{N}_+}$ by
\[
s_1 = \frac{1}{a_1}, \quad s_{n+1} = \frac{1}{s_n + a_{n+1}}, \ n \in \mathbb{N}_+.
\]
Let us consider the RSCC $\left\{(W, {\mathcal W}), (X, {\mathcal X}), u, P\right\}$, where 
\[
W = [0,1], \quad {\mathcal W} = {\mathcal B}_{[0,1]},
\]
\[
X = \mathbb{N}_+, \quad {\mathcal X} = {\mathcal P}_{\mathbb{N}_+}, 
\]
\[
u: W \times X \rightarrow W, \quad u(w, x) = \frac{1}{w+x}, 
\]
\[
P : W \times {\mathcal X} \rightarrow W, \quad P(w, x) = \frac{w+1}{(w+x)(w+x+1)}.
\]
The sequences $(a_n)_{n \in \mathbb{N}_+}$ and $(s_n)_{n \in \mathbb{N}_+}$, $s_0 = 0$, are equivalent to the chain with complete connections $(\zeta_n)_{n \in \mathbb{N}_+}$ and the Markov chain $(\xi_n)_{n \in \mathbb{N}}$ associated with the above RSCC. More precisely, defining the one-to-one map $\theta$ from $\left(\mathbb{N}_+\right)^{\mathbb{N}_+}$ into $[0,1]$ by 
\[
\theta(a_1, a_2, a_3, \ldots) = \displaystyle \frac{1}{a_1+\displaystyle \frac{1}{a_2+\displaystyle \frac{1}{a_3 + \ddots}}}, \quad a_i \in \mathbb{N}_+, \ i \in \mathbb{N}_+,
\]
we have $\zeta_n(\sigma) = a_n(\theta(\sigma))$, $\xi_n(\sigma) = s_n(\theta(\sigma))$, $n \in \mathbb{N}_+$, $\sigma \in \left(\mathbb{N}_+\right)^{\mathbb{N}_+}$.
\hfill $\Box$
\end{example} 

\subsubsection{Properties of the associated operators}

In this subsection we present the asymptotic and ergodic properties of the associated operators. These properties are used to obtain the ergodicity of a RSCC by letting the associated Markov chain satisfy some topological properties. To state these results we need some preliminary definitions. 

Let $Q_n$ be the transition probability function defined by 
\[
Q_n(w, A) = \frac{1}{n}\sum_{k=1}^{n} Q^k(w, A)
\]
where $Q^k$, $k \geq 1$, is the $k$-step transition probability function of the Markov chain associated with RSCC (\ref{5.1}). Let $U_n$ be the Markov operator associated with $Q_n$. 

Next, let us consider the norm $\left\| \cdot \right\|_L$ defined on $L(W)=$ the space of Lipschitz complex-valued functions defined on $W$ by
\[
\left\| f \right\|_L = \sup_{w \in W} |f(w)| + \sup_{w'\ne w''} \frac{|f(w') - f(w'')|}{|w' - w''|}, \quad f \in L(W).
\]
As is well known, $(L(W), \left\| \cdot \right\|_L)$ is a Banach space.

The following can be found in \cite{r10}.

If there exists a linear bounded operator $U^\infty$ from $L(W)$ to $L(W)$ such that
\[
\lim_{n \rightarrow \infty} \left\|U_nf-U^\infty f\right\|_L = 0, 
\]
for any $f \in L(W)$ with $\left\|f\right\|_L=1$, we say $U$ \textit{ordered}.

If
\[
\lim_{n \rightarrow \infty} \left\|U^nf-U^\infty f\right\|_L = 0,
\]
for any $f \in L(W)$ with $\left\|f\right\|_L=1$, we say $U$ \textit{aperiodic}, where $U^n$ is the $n$th iterate of $U$, $n \in \mathbb{N}$, with $U^0$ is the identity. 

If $U$ is ordered and $U^\infty(L(W))$ is one-dimensional space, it is named \textit{ergodic} with respect to $L(W)$.

If $U$ is ergodic and aperiodic, it is named \textit{regular} with respect to $L(W)$ and the corresponding Markov chain has the same name.

The definition below is due to M.F. Norman \cite{N-1972} and isolates a class of RSCCs, called \textit{RSCCs with contraction}.

An RSCC $\left\{(W, {\mathcal W}), (X, {\mathcal X}), u, P\right\}$  is said to be an RSCC \textit{with contraction} if and only if there is a distance $\mathrm{d}$ on $W$ and the metric space $(W, \mathrm{d})$ is separable, $r_1< \infty$, $R_1< \infty$, and there exists a natural integer $k$ such that $r_k<1$, where 
\[
r_k = \sup_{w'\neq w''} \sum_{X^k}P_k\left(w, x^{(k)}\right)\frac{\mathrm{d}\left(w'x^{(k)}, w''x^{(k)}\right)}{\mathrm{d}\left(w', w''\right)}, \quad k \in \mathbb{N}_+, 
\]
and
\[
R_k = \sup_{A \in {\mathcal X}^k} \sup_{w'\neq w''} \frac{P_k\left(w', A\right)- P_k\left(w'', A\right)}{d\left(w', w''\right)}.
\]
The following result can be found in \cite{r10}. 

\begin{theorem} \label{th.5.2}
Let $W$ be a compact metric space with a distance $\mathrm{d}$ and $\left\{(W, {\mathcal W}), (X, {\mathcal X}), u, P\right\}$ be a RSCC with contraction. 
\begin{enumerate}
\item[(i)] 
The Markov chain associated to the RSCC is regular if and only if there exists a point $w_0 \in W$ such that
\[
\lim_{n \rightarrow \infty} \mathrm{d}\left(\sigma_n(w), w_0\right)=0, 
\]
for any $w \in W$, where $\sigma_n(w) = \mathrm{supp \ } Q^n(w, \cdot)$ ($\mathrm{supp \ }\mu$ denotes the support of the measure $\mu$). 

\item[(ii)]
The suports of $Q^n(w, \cdot)$, $n \in \mathbb{N}_+$, $w \in W$, can be iteratively computed as follows: 
\[
\sigma_{m+n}(w) = \overline{\bigcup_{w' \in \sigma_m(w)}\sigma_n(w')}, 
\]
for any $m$, $n \in \mathbb{N}_+$, $w \in W$, where the overline means the topological closure.
\end{enumerate}
\end{theorem}

An RSCC $\left\{(W, {\mathcal W}), (X, {\mathcal X}), u, P\right\}$, whose associated Markov chain is regular with respect to $B((W, {\mathcal W}))$, is \textit{uniformly ergodic} and $\displaystyle\lim_{n \rightarrow \infty}\varepsilon_n = 0$, where
\[
\varepsilon_n := \sup_{\begin{array}{cc}
	w \in W, r \in \mathbb{N}_+\\
	A \in {\mathcal X}^r\\ 
\end{array}} \left|P^n_r(w,A) - \mathrm{P}^\infty_r(A)\right|,
\]
while \rmfamily{P}$^\infty_r$ is the probability on ${\mathcal X}^r$.

\begin{theorem} \label{th.5.4}
Let $W$ be a compact metric space with a distance $\mathrm{d}$. If the RSCC $\left\{(W, {\mathcal W}), (X, {\mathcal X}), u, P\right\}$ with contraction has regular associated Markov chain, then it is uniformly ergodic.
\end{theorem}

\subsection{The RSCC associated with expansion of the type of (\ref{3.1})}

First, it is easy to check that $P_m^i$ from $(\ref{3.40})$ defines a transition probability function from $\left(I,{\mathcal{B}}_I\right)$ to $\left(\mathbb{N}, {\mathcal{P}}(\mathbb{N}) \right)$, i.e., $\displaystyle\sum_{i \in \mathbb{N}}P_m^i(x)=1$, $x \in I$.

Let us to consider the random system with complete connections
\begin{equation}
\left\{\left(I, {\mathcal{B}}_I\right), \left(\mathbb{N}_+, {\mathcal{P}}(\mathbb{N}_+)\right), u, P\right\}, \label{72}
\end{equation}
where $u:I \times \mathbb{N} \rightarrow I$, $u(x,i)=u^i_m(x)$ is given in (\ref{3.63}) and the function $P(x,i) = P^i_m (x)$ given in (\ref{3.40}). 

We denote by $U_m$ the associated Markov operator of RSCC (\ref{72}) with the transition probability function 
\[
Q_m(x, A) = \sum_{\left\{i \in \mathbb{N}: u^i_m(x) \in A \right\}}P^i_m(x), \quad x \in I,  \mbox{ } A \in {\mathcal{B}}_I. 
\]
Then $Q_m^n(\cdot, \cdot)$ will denote the $n$-step transition probability function of the same Markov chain.

The ergodic behaviour of RSCC (\ref{72}) allows us to find the limiting distribution function $F$ and the invariant measure $\mathrm{Q}_m^{\infty}$ induced by $F$. 

\begin{proposition} \label{3.5.2} 
RSCC $(\ref{72})$ is uniformly ergodic.
\end{proposition} 
\textit{Proof.} We apply Theorem \ref{th.5.4}. Putting $\Delta_i = m^{-i} - m^{-2i}$, $i \in \mathbb{N}$, we get
\[
P^i_m (x) = (m-1) \left[m^{-(i+1)} + \frac{\Delta_i}{x+(m-1)m^{-i}+1} - \frac{\Delta_{i+1}}{x + (m-1)m^{-(i+1)}+1}\right], 
\]
We have 
\begin{eqnarray*}
\frac{d}{dx}u(x,i) &=& - \frac{(m-1)m^{-i}}{((m-1)x+1)^2}, \\
\frac{d}{dx}P(x,i) &=& (m-1) \left[\frac{\Delta_{i+1}}{(x + (m-1)m^{-(i+1)}+1)^2} - \frac{\Delta_i}{(x+(m-1)m^{-i}+1)^2} \right],
\end{eqnarray*}
for all $x \in I$ and $i \in \mathbb{N}$, so that $\displaystyle\sup_{x \in I}\left|\frac{d}{dx}u(x,i)\right| = (m-1)m^{-i}$ and $\displaystyle\sup_{x \in I}\left|\frac{d}{dx}P(x,i)\right| < \infty$. 
Hence the requirements of definition of an RSCC with contraction are fulfilled. To prove the regularity of $U$ with respect to $L(I)$ let us define recursively $x_{n+1} = (x_n + 2)^{-1}$, $n \in \mathbb{N}$, with $x_0 = x$.

A criterion of regularity is expressed in Theorem \ref{th.5.2}(i), in terms of supports $\sigma_n(x)$ of the $n$-step transition probability functions $Q_m^n(x, \cdot)$, $n \in \mathbb{N}_+$. Clearly $x_{n+1} \in \sigma_1(x_n)$ and therefore Theorem \ref{th.5.2}(ii) and an induction argument lead us to the conclusion that $x_n \in \sigma_n(x)$, $n \in \mathbb{N}_+$. But, $\displaystyle \lim_{n \rightarrow \infty} x_n = \sqrt{2}-1$ for any $x \in I$. Hence
\[
\mathrm{d}\left( \sigma_n(x), \sqrt{2}-1 \right) \leq \left|x_n - \sqrt{2}+1\right| \rightarrow 0 \ \mbox{ as } \ n \rightarrow \infty,
\]
where $\mathrm{d}(x, y) = \left|x - y\right|$, for any $x, y \in I$. The regularity of $U_m$ with respect to $L(I)$ follows from Theorem \ref{th.5.2}. Moreover, $Q_m^n (\cdot, \cdot)$ converges uniformly to a probability measure $Q_m^\infty$ and that there exist two positive constants $q < 1$ and $k$ such that
\begin{equation}
\left\|U^n_mf-U^\infty_m f\right\|_L \leq kq^n\left\|f\right\|_L, \quad n \in \mathbb{N}_+, \mbox{ } f \in L(I), \label{74}
\end{equation}
where
\begin{eqnarray}
U^n_mf(\cdot) &=& \int_{I}f(y)Q^n_m(\cdot, dy), \\
U^\infty_m f &=& \int_{I}f(y)\mathrm{Q}^\infty_m(dy), \label{75}
\end{eqnarray}
and $\mathrm{Q}^\infty_m$ is the invariant probability measure of the transformation $\tau_m$, i.e., $\mathrm{Q}^\infty_m$ has the the density $\rho_m (x)$ given in $(\ref{3.41})$, $x \in I$. 

\hfill $\Box$
 
Now we are able to find the limiting distribution function
\[
F(x) = F_\infty (x) = \lim_{n \rightarrow \infty} \mu (\tau_m^n < x)
\]
and obtain a convergence rate result.

\subsection{Proof of Theorem \ref{G-K}}
We prove Theorem \ref{G-K} in this subsection.

\noindent \textbf{Proof of Theorem \ref{G-K}}

By (\ref{75}) we have
\[
U^\infty_m f_0 = \int_{I} f_0(y)\mathrm{Q}^\infty_m(dx) = k_m, \quad f_0 \in L(I). 
\]
Taking into account (\ref{74}), there exist two constants $q<1$ and $k$ such that
\[
\left\|U^n_mf_0 - U^\infty_m f_0\right\|_L \leq k q^n\left\|f_0\right\|_L, \quad n \in \mathbb{N}_+.
\]
Further, consider $C(I)$ the metric space of real continuous functions defined on $I$ with the supremum norm $\displaystyle \left\|f\right\| = \sup_{x \in I}\left|f(x)\right|$. Since $L(I)$ is a dense subset of $C(I)$ we have
\begin{equation}
\lim_{n \rightarrow \infty} \left\|\left(U^n_m - U^\infty_m\right)f_0\right\| = 0, \label{90}
\end{equation}
for all $f_0 \in C(I)$. Therefore, (\ref{90}) is valid for a measurable function $f_0$ which is $\mathrm{Q}^\infty_m$-almost surely continuous, that is, for a Riemann-integrable function $f_0$. Thus, we have
\begin{eqnarray*}
F(x) &=& \lim_{n \rightarrow \infty} \mu \left(\tau^n_m < x\right) = \lim_{n \rightarrow \infty} \int_{0}^{x} U^n_mf_0(u) \rho_m(u) du  \\
&=& k_m \int_{0}^{x} \rho_m(u)du \\
&=& \frac{k_m}{(m-1)^2}\log \frac{m((m-1)x+1)}{(m-1)x+m}. 
\end{eqnarray*}
Hence \ref{86} is proved.
\hfill $\Box$
\\

\noindent \textbf{Acknowledgments} 

The author would like to express their sincere thanks to the referees for their valuable comments.

\appendix
\section{Proofs of propositions}

We prove Propositions \ref{prop.3.3.1}, \ref{prop.3.1}, \ref{prop.4.1} and \ref{prop.4.2} in this section.
\\

\noindent \textbf{Proof of Proposition \ref{prop.3.3.1}}

We briefly give some general properties about the Perron-Frobenius operator (see, e.g., \cite{r1, r11}) which will be useful both to demonstrate this proposition and in Section 4.

Let $(X, \mathcal{X}, \mu)$ be a probability space. A transformation $\tau$ of $X$ is said to be $\mu$\textit{-non-singular} if and only if $\mu \left(\tau^{-1}(A)\right) = 0$ for all $A \in \mathcal{X}$ for which $\mu(A) = 0$; it is said to be \textit{measure-presearving} if and only if $\mu \tau^{-1} = \mu$, i.e., $\mu\tau^{-1}(A) = \mu(A)$ for all $A \in \mathcal{X}$. Clearly, any $\mu$-preserving transformation is $\mu$- non-singular.

The Perron-Frobenius operator $P_{\mu}$ associated with a $\mu$-non-singular transformation $\tau$ is defined as the linear bounded operator on $L^1_{\mu} = \left\{f:I \rightarrow \mathbb{C} : \int_{I}\left|f\right|d\mu<\infty\right\}$ which takes $f \in L^1_{\mu}$ into $P_{\mu}f \in L^1_{\mu}$ with 
\[
\int_{A}P_{\mu}f \mbox{d}\mu = \int_{\tau^{-1}(A)}f \mbox{d}\mu, \quad A \in {\mathcal{X}},
\]
or, equivalently
\[
\int_{X}gP_{\mu}f \mbox{d}\mu = \int_{X} (g \circ \tau)f \mbox{d}\mu
\]
for all $f \in L^1_{\mu}$ and $g \in L^{\infty}_{\mu}$. 

In particular, the Perron-Frobenius operator $P_{\lambda}$ of $\tau$ under the Lebesgue measure $\lambda$ is (see \cite{r1}, p.86) 
\begin{equation}
P_{\lambda}f(x)=\frac{d}{dx}\int_{\tau^{-1}([0,x])}f d\lambda=\sum_{t \in \tau^{-1}(x)}\frac{f(t)}{\left|\tau'(t)\right|} \mbox{  } \mbox{  } \mbox{a.e. in } I. \label{3.61}
\end{equation}
The probabilistic interpretation of $P_{\mu}$ is immediate: if an $X$-valued random variable $\xi $ on $X$ has $\mu$-density $h$, that is, $\mu(\xi \in A) = \int_{A}h \mbox{d}\mu$, $A \in X$, with $h \geq 0$ and $\int_{X}h \mbox{d}\mu=1$, then $\tau \circ \xi$ has $\mu$-density $P_{\mu}h$. The following properties hold:
\begin{enumerate}
\item[(i)] $P_{\mu}$ is positive, that is, $P_{\mu}f \geq 0$ if $f \geq 0$;
\item[(ii)] $P_{\mu}$ preserves integrals, that is, $\int_{X} P_{\mu}f\mbox{d}\mu = \int_{X} f \mbox{d}\mu$, $f \in L^1_{\mu}$;
\item[(iii)] $\left\|P_{\mu}\right\|_{p,\mu} := \sup \left(\left\|P_{\mu}f\right\|_{p,\mu}: f \in L^p_{\mu}, \left\|f\right\|_{p,\mu} = 1\right) \leq 1$ for any $p \geq 1$ and $p = \infty$;
\item[(iv)] for any $n \in \mathbb{N}_+$ the $n$th power $P_{\mu}^n$ of $P_{\mu}$ is the Perron-Frobenius operator associated with the $n$th iterate $\tau^n$ of $\tau$ under $\mu$;
\item[(v)] $(P_{\mu}f)^* = P_{\mu}f^*$ for any $f \in L^1_{\mu}$, where $z^*=$ complex conjugate of $z \in \mathbb{C}$ (=the set of complex numbers);
\item[(vi)] $P_{\mu}((g \circ \tau)f)= gP_{\mu}f$ for any $f \in L^1_{\mu}$ and $g \in L^{\infty}_{\mu}$;
\item[(vii)] $P_{\mu}f = f$ if and only if $\tau$ is $\nu$-preserving, where $\nu$ is defined by $\nu(A) = \int_{A} f\mbox{d}\mu$, $A \in \mathcal{X}$. In particular, $P_{\mu}1=1$ if and only if $\tau$ is $\mu$-preserving.
\end{enumerate}

\textit{Proof of the Proposition \ref{prop.3.3.1}} From above, it is sufficient to show that the function $\rho_m$ defined in (\ref{3.41}) is an eigenfunction of the Perron-Frobenius operator of $\tau_m$ with the eigenvalue $1$: 
\begin{equation}
P_{\tau_m}\rho_m(x) = \sum_{t \in \tau^{-1}_m(x)}\frac{\rho_m(t)}{\left|\tau'_m(t)\right|}. \label{3.42}
\end{equation}
First, we note that
\begin{equation}
\tau^{-1}_m(x) = \left\{ \frac{m^{-i}}{1+(m-1)x}: i \geq 1, x \in I  \right\}. \label{3.43}
\end{equation}
Thus
\begin{eqnarray*}
P_{\tau_m}\rho_m(x) = \sum_{i=0}^{\infty}\frac{(m-1)m^{-i}}{(1+(m-1)x)^2}\rho_m\left(\frac{m^{-i}}{1+(m-1)x}\right) \qquad \qquad \qquad \qquad \\
= \sum_{i=0}^{\infty} (m-1)m^{-(i+1)} \frac{1}{\left((m-1)x+(m-1)m^{-(i+1)}+1\right)} \qquad \ \ \ \\ 
\times \frac{1}{\left((m-1)x+(m-1)m^{-i}+1\right)} \qquad \qquad \qquad \qquad \qquad \quad \\
= \frac{1}{m-1} \sum_{i=0}^{\infty} \left(\frac{1}{(m-1)x+(m-1)m^{-(i+1)}+1} \right. \qquad \qquad \quad \ \\
- \left. \frac{1}{(m-1)x+(m-1)m^{-i}+1}\right) \qquad \qquad \qquad \qquad \qquad \ \ \ \ \\
= \frac{1}{m-1}\left(\frac{1}{(m-1)x+1} - \frac{1}{(m-1)x+m} \right) \qquad \qquad \qquad \ \ \ \\
= \frac{1}{((m-1)x+1)((m-1)x+m)} = \rho_m(x). \qquad \qquad \qquad \ \ \ 
\end{eqnarray*}

\hfill $\Box$
\\

\noindent \textbf{Proof of Proposition \ref{prop.3.1}}

We should show that $\overline{\gamma}_m \left(\overline{\tau}^{-1}_m(B)\right) = \overline{\gamma}_m (B)$ for any $B \in {\mathcal{B}}^2_I$ or, equivalently, since $\overline{\tau}_m$ is invertible on $\Omega^2$, that
\begin{equation}
\overline{\gamma}_m \left(\overline{\tau}_m(B)\right) = \overline{\gamma}_m (B), \mbox{ for any } B \in {\mathcal{B}}^2_I. \label{3.51}
\end{equation}
We start with $B = (a, b) \times (c, d)$, where 
\[
a = m^{-(i+1)}, \mbox{ } b = m^{-i}, \quad i \in \mathbb{N}
\]
and $c$ and $d$ arbitrary numbers from $(0, 1)$. Then
\begin{equation}
\overline{\tau}_m(B) = \left\{\left( \tau_m(x), \frac{m^{-a_1(x)}}{(m-1)y+1} \right) \left| x \in (a, b), y \in (c, d)\right.\right\}. \label{3.52}
\end{equation}
Taking $x = m^{-(i+\theta)}$, $0<\theta<1$, we have
\[
\tau_m(x) = \frac{m^{\theta} - 1}{m - 1}, \quad a_1(x) = i
\]
such that
\begin{equation}
\overline{\tau}_m(B) = \left((0, 1), \left( \frac{m^{-i}}{(m-1)d+1}, \frac{m^{-i}}{(m-1)c+1} \right)\right). \label{3.53}
\end{equation}
Let
\[
I(m,i,c,d) \equiv \left( \frac{m^{-i}}{(m-1)d+1}, \frac{m^{-i}}{(m-1)c+1} \right).
\]
A simple computation yields
\begin{eqnarray}
\overline{\gamma}_m \left(\overline{\tau}_m(B)\right) &=& k_m \int^{1}_{0}dx \int_{I(m,i,c,d)}\frac{dy}{((m-1)(x+y)+1)^2} \nonumber \\
&=& k_m \int^{m^{-i}}_{m^{-(i+1)}}dx \int^{d}_{c} \frac{dy}{((m-1)(x+y)+1)^2} = \overline{\gamma}_m (B) \nonumber
\end{eqnarray}
that is, (\ref{3.51}) holds. 

Next, we consider the case 
\[
a = \frac{m^{-i}}{(m-1)m^{-j}+1}, \quad b = \frac{m^{-i}}{(m-1)m^{-(j+1)}+1}, \quad i, \mbox{ } j \in \mathbb{N}
\]
and $(c, d)$ an arbitrary interval. Now, with 
\[
x = \frac{m^{-i}}{(m-1)m^{-(j+\theta)}+1},
\]
we have
\[
\left\{\frac{\log x^{-1}}{\log m}\right\} = \left\{i + \frac{\log \left(1+(m-1)m^{-(j+\theta)}\right)}{\log m}  \right\} = \frac{\log \left(1+(m-1)m^{-(j+\theta)}\right)}{\log m} 
\]
and
\[
a_1(x) = \left\lfloor \frac{\log x^{-1}}{\log m} \right\rfloor = i.
\]
Thus, 
\[
(m-1)\tau_m(x) = m^{\frac{\log \left(1+(m-1)m^{-(j+\theta)}\right)}{\log m}} - 1 = (m-1) m^{-(j+\theta)}.
\]
Hence, 
\begin{equation}
\overline{\tau}_m(B) = \left(m^{-(j+1)}, m^{-j}\right) \times \left( \frac{m^{-i}}{(m-1)d+1}, \frac{m^{-i}}{(m-1)c+1} \right). \label{3.54}
\end{equation}
A straightforward calculation shows us that
\begin{eqnarray}
\overline{\gamma}_m \left(\overline{\tau}_m(B)\right) &=& k_m \int^{m^{-j}}_{m^{-(j+1)}}dx \int_{I(m,i,c,d)} \frac{dy}{((m-1)(x+y)+1)^2} \nonumber \\
&=& k_m \int_{I(m,i,m^{-j},m^{-(j+1)})}dx \int^{d}_{c} \frac{dy}{((m-1)(x+y)+1)^2} = \overline{\gamma}_m (B) \nonumber
\end{eqnarray}
that is, (\ref{3.51}) holds.

Since any arbitrary interval $(a, b)$ can be written as a reunion of fundamental intervals the proof is complete.
\hfill $\Box$
\\

\noindent \textbf{Proof of Proposition \ref{prop.4.1}}
(i) Let $\tau_{m,i} : I_i \rightarrow I$ denote the restriction of $\tau_m$ to the interval $I_i=\left(m^{-(i+1)}, m^{-i}\right]$, $i \in \mathbb{N}$, that is, 
\begin{equation}
\tau_{m,i}(x) = \frac{1}{m-1}\left(\frac{m^{-i}}{x} - 1\right), \mbox{ } x \in I_i. \label{3.64}
\end{equation}
For any $f \in L^1_{\gamma_m}$ and any $A \in {\mathcal{B}}_I$, we have
\begin{equation}
\int_{\tau^{-1}_m(A)} f d \gamma_m = \sum_{i \in \mathbb{N}} \int_{\tau^{-1}_m(A \cap I_i)}f d \gamma_m = \sum_{i \in \mathbb{N}} \int_{\tau^{-1}_{m,i}(A)}f d \gamma_m. \label{3.65}
\end{equation}
For any $i \in \mathbb{N}$, by the change of variable
\begin{equation}
x = \tau^{-1}_{m,i}(y) = \frac{m^{-i}}{(m-1)y+1}, \label{3.66}
\end{equation}
we successively obtain
\begin{eqnarray}
\int_{\tau^{-1}_{m,i}(A)}f d \gamma_m &=& k_m \int_{\tau^{-1}_{m,i}(A)} \frac{f(x)}{((m-1)x+1)((m-1)x+m)}dx \nonumber \\
&=& k_m \int_{A} \frac{f\left(u^i_m(y)\right)}{\left((m-1)u^i_m(y)+1\right)\left((m-1)u^i_m(y)+m\right)} \nonumber \\
&\times& \frac{(m-1)m^{-i}}{((m-1)y+1)^2}dy \nonumber \\
&=& k_m \int_{A} f\left(u^i_m(y)\right)(m-1)m^{-(i+1)} \frac{1}{\left((m-1)y+(m-1)m^{-i}+1\right)} \nonumber \\
&\times& \frac{1}{\left((m-1)y+(m-1)m^{-(i+1)}+1\right)} dy \nonumber \\
&=& \int_{A} P^i_m((m-1)y)f\left(u^i_m(y)\right)\gamma_m (dy). \label{3.67}
\end{eqnarray}
Now, (\ref{3.62}) follows from (\ref{3.65}) and (\ref{3.67}).
\hfill $\Box$

(ii)(a) From (\ref{3.64}) and (\ref{3.66}), for any $f \in L^1_{\gamma_m}$ and any $A \in {\mathcal{B}}_I$, we have
\[
\begin{array}{l}
\displaystyle \int_{\tau^{-1}_m(A)} f \mbox{d} \mu =\displaystyle  \sum_{i \in \mathbb{N}} \displaystyle \int_{\tau^{-1}_m(A \cap I_i)}f \mbox{d}\mu = \displaystyle \sum_{i \in \mathbb{N}} \displaystyle \int_{\tau^{-1}_{m,i}(A)}f \mbox{d}\mu \nonumber \\
= \displaystyle \sum_{i \in \textbf{N}}\displaystyle \int_{\tau^{-1}_{m,i}(A)}f(x)h(x)\mbox{d}x = \displaystyle \sum_{i \in \textbf{N}} \displaystyle \int_{A} \displaystyle \frac{f(u_m^i(y))h(u_m^i(y))(m-1)m^{-i}}{((m-1)y+1)^2}\mbox{d}y \nonumber \\
\end{array}
\]
\begin{equation}
\qquad \qquad \ \ \ = \int_{A} \sum_{i \in \textbf{N}} \frac{h(u_m^i(x))}{((m-1)x+1)^2}(m-1)m^{-i}f(u_m^i(x))\mbox{d}x. \label{3.71} 
\end{equation}
Since $\mbox{d}\mu = h \mbox{d}\lambda$, (\ref{3.68}) follows from (\ref{3.71}). Now, since $g(x)=((m-1)x+1)((m-1)x+m)f(x)h(x)$, from (\ref{3.62}) we have 
\begin{equation}
U_mg(x) = \frac{((m-1)x+m)}{(m-1)x+1}(m-1) \sum_{i \in \mathbb{N}}m^{-i}h(u_m^i(x))f(u_m^i(x)). \label{3.72}
\end{equation}
Now, (\ref{3.69}) follows immediately from (\ref{3.68}) and (\ref{3.72}).
\hfill $\Box$

(b)We will use mathematical induction. For $n=0$, the equation (\ref{3.75}) reduces to
\begin{equation}
\mu (A) = \int_{A} h(x)dx, \mbox{ } A \in {\mathcal{B}}_I, \nonumber
\end{equation}
which is obviously true. Assume that (\ref{3.75}) holds for some $n \in \mathbb{N}$. Then 
\begin{eqnarray}
\mu \left(\tau^{-(n+1)}_m(A)\right) &=& \mu \left(\tau^{-n}_m(\tau^{-1}_m(A))\right) \nonumber \\
&=& \int_{\tau^{-1}_m(A)} \frac{U^n_mf(x)}{((m-1)x+1)((m-1)x+m)}dx \nonumber \\
&=& \frac{1}{k_m}\int_{\tau^{-1}_m(A)} U^n_mf(x) d \gamma_m(x). \nonumber
\end{eqnarray}
By the very definition of the Perron-Frobenius operator $U_m=P_{\gamma_m}$ we have
\begin{equation}
\int_{\tau^{-1}_m(A)} U^n_m f d \gamma_m = \int_{A} U^{n+1}_m f d \gamma_m. \nonumber
\end{equation}
Therefore, 
\begin{eqnarray}
\mu \left(\tau^{-(n+1)}_m(A)\right) &=& \frac{1}{k_m}\int_{A} U^{n+1}_mf d \gamma_m \nonumber \\ 
&=& \int_{A} \frac {U^{n+1}_mf(x)}{((m-1)x+1)((m-1)x+m)}dx \nonumber
\end{eqnarray}
which ends the proof.
\hfill $\Box$
\\
\noindent \textbf{Proof of Proposition \ref{prop.4.2}}
(i) For $x, y \in I$ we have
\begin{eqnarray*}
U_mf(x) - U_mf(y) &=& \sum_{i \in \mathbb{N}} (P_m^i((m-1)x)f(u_m^i(x)) - P_m^i((m-1)y)f(u_m^i(y))) \\
                  &=& \sum_{i \in \mathbb{N}} (P_m^i((m-1)x - P_m^i((m-1)y)(f(u_m^i(x)) - f(u_m^0(x))) \\
                  &+& \sum_{i \in \mathbb{N}} P_m^i((m-1)y (f(u_m^i(x)) - f(u_m^i(y))) \\
                  &=& \sum_{i \in \mathbb{N}_+}(P_m^i((m-1)x - P_m^i((m-1)y)(f(u_m^i(x)) - f(u_m^0(x))) \\
                  &+& \sum_{i \in \mathbb{N}} P_m^i((m-1)y) (f(u_m^i(x)) - f(u_m^i(y))).  
\end{eqnarray*}
Note that the function $P_m^0$ is increasing, while the functions $P_m^i$, $i \in \mathbb{N}_+$, are all decreasing. Let $x < y$, with $x, y \in I$. It follows from the above equation that 
\begin{eqnarray*}
\left|U_mf(x) - U_mf(y)\right| &\leq& \left( \sum_{i \in \mathbb{N}_+}(P_m^i((m-1)x - P_m^i((m-1)y)\right)\mbox{var } f \\
                               &+& \sup_{y \in I, i \in \mathbb{N}} P_m^i((m-1)y) \sum_{i \in \mathbb{N}} \mbox{var}_{[x,y]} f \circ u_i(x) \\
                               &=& (1 - P_m^0((m-1)x) - 1 + P_m^0((m-1)y)) \mbox{var } f \\
                               &+& P_m^0(m-1) \sum_{i \in \mathbb{N}} \mbox{var}_{[x,y]} f \circ u_i(x).
\end{eqnarray*}
Hence
\begin{eqnarray*}
\mbox{var } U_mf \leq (2P_m^0(m-1) - P_m^0(0)) \mbox{var} f = \left( \frac{2m(m-1)}{m^2+m-1} - \frac{m-1}{2m-1}\right) \mbox{var } f \\
                = \frac{(m-1)\left(3m^2-3m+1\right)}{(2m-1)\left(m^2+m-1\right)} \mbox{var } f. \ \ 
\end{eqnarray*}
Define $f$ by $f(x)=0$, $0 \leq x \leq \frac{1}{m}$, and $f(x)=1$, $\frac{1}{m} < x \leq 1$. Then we have $U_mf(x) = P_m^0(x)$, $0 \leq x < 1$ and $U_mf(1) = 0$. Since $\mbox{var } U_mf = \displaystyle\frac{(m-1)\left(3m^2-3m+1\right)}{(2m-1)\left(m^2+m-1\right)}$ and $\mbox{var } f = 1$, it follows that the constant $K_m$ cannot be lowered.
\hfill $\Box$

(ii) The transition operator of $(s^a_n)_{n \in \mathbb{N}_+}$ takes $f \in B(I)$ to the function defined by 
\begin{eqnarray}
E_a\left( \left. f(s^a_{n+1})\right| s^a_n = s \right) &=& \sum_{i \in \mathbb{N}}P_m^i((m-1)s)f(u_m^i(s)) \nonumber \\
&=& U_m f(s), \quad s \in I,
\end{eqnarray}
where $E_a$ stands for the mean value operator with respect to the probability measure $\gamma^a_m$.

\hfill $\Box$



\begin{thebibliography}{[01]}

\bibitem{BE-1988} Barnsley, M., Elton, J., \textit{A New Class of Markov Processes for Image Encoding}, Adv. in Appl. Probab. \textbf{20} (1988), 14–32.
\bibitem{Barnsley-1988} Barnsley, M., Demko, S., Elton, J., Gerinomo, J., \textit{Invariant measures for Markov processes arising from iterated function systems with place-dependent probabilities}, Ann. Inst. H. Poincar\'e, Probab. Statist \textbf{24}(3) (1988), 367-394.
\bibitem{BK-1990} Bosma, W., Kraaikamp, C., \textit{Metrical Theory for Optimal Continued Fractions}, J. Number Theory \textbf{34} (1990), 251-270.
\bibitem{r1} Boyarsky, A., G\'ora, P., \textit{Laws of Chaos: Invariant Measures and Dynamical Systems in One Dimension.} Birkh\"auser, Boston, 1997.
\bibitem{r5} Chan, H.-C., {\it The asymptotic growth rate of random Fibonacci type sequences. II}, Fibonacci Quart.\textbf{44} (2006), 73-84. 
\bibitem{r7} Doeblin, W., Fortet, R., \textit{Sur des cha\^ines \`a liaisons compl\`etes}, Bull. Soc. Math. France \textbf{65} (1937), 132-148.
\bibitem{r8} Harris, T.E., \textit{On chains of infinite order.} Pacific J. Math. \textbf{5} (1955), 707-724.
\bibitem{HIR-2003} Herkenrath, U., Iosifescu, M., Rudolph, A., \textit{Random systems with complete connections and iterated function systems}, Math. Rep. \textbf{5}(55) (2003), 127–140.
\bibitem{MI-1963}Iosifescu, M., {\it Random systems with complete connections with an arbitrary set of states}, 
Rev. Roumaine Math. Pures Appl. \textrm{8} (1963), 611-645.
\bibitem{r10} Iosifescu, M., Grigorescu, S., \textit{Dependence With Complete Connections and its Applications.} Cambridge Tracts in Mathematics \textbf{96}, 1990. Cambridge Univ.Press, Cambridge. [(2009): second printing slightly corrected].
\bibitem{r11} Iosifescu, M., Kraaikamp, C., {\it Metrical Theory of Continued Fractions.} Kluwer Academic Publisher, Dordrecht, 2002.
\bibitem{r13} Iosifescu, M., Theodorescu, R., \textit{Random Processes and Learning.} Springer-Verlag, Berlin, 1969.
\bibitem{r14} Kalpazidou, S., {\it On a problem of Gauss-Kuzmin type for continued fraction with odd partial quotients}, Pacific J. Math. \textbf{123} (1) (1986), 103-114.
\bibitem{Karlin-1953} Karlin, S., \textit{Some random walks arising in learning models. I}, Pacific J. Math. \textbf{3} (4) (1953), 725-756.
\bibitem{Hincin-1964} Khinchin, A.Ya., \textit{Continued Fractions.} Univ. Chicago Press, Chicago, 1964 [Translation of the 3rd (1961) Russian Edition].
\bibitem{r17} Kuzmin, R.O., \textit{On a problem of Gauss}, Dokl. Akad. Nauk SSSR Ser. A (1928) 375-380. [Russian; French version in \textit{Atti Congr. Internaz.Mat.} (\textit{Bologna}, 1928), Tomo \textbf{VI} (1932) 83-89. Zanichelli, Bologna]. 
\bibitem{r19} L\'evy, P., \textit{Sur les lois de probabilit\'e dont d\'ependent les quotients complets et incomplets d'une fraction continue}, Bull. Soc. Math. France \textbf{57} (1929), 178-194. 
\bibitem{MU-2003} Mauldin, R.D., Urba\'nski, M., \textit{The doubling property of conformal measures of infinite iterated function systems}, J. Number Theory \textbf{102}  (1) (2003), 23-40.
\bibitem{r20} Nakada, H., \textit{Metrical theory for a class of continued fraction transformations and their natural extensions}, Tokyo J.Math. \textbf{4}(2) (1981), 399-426. 
\bibitem{N-1972} Norman, E., \textit{Markov Processes and Learning Models.} Academic Press, New York, 1972.
\bibitem{r21} Onicescu, O., Mihoc, Gh., \textit{Sur les cha\^ines de variables statistiques.} Bull. Sci.Math. \textbf{59} (1935), 174 - 192.
\bibitem{Sebe-2001} Sebe, G.I., {\it On convergence rate in the Gauss-Kuzmin problem for the grotesque continued fractions},
Monatsh. Math. \textbf{133} (3) (2001), 241-254.
\bibitem{Sebe-2002} Sebe, G.I., {\it A Gauss-Kuzmin theorem for the Rosen fractions}, J. Th\'eor. Nombres Bordeaux \textbf{14} (2) (2002), 667-682.
\bibitem{r23} Sebe, G.I., {\it On a Gauss-Kuzmin-type problem for a new continued fraction expansion with explicit invariant measure}, Proc. of the 3-rd Int. Coll "Math. in Engg.and Numerical Physics" (MENP-3), 7-9 October 2004 Bucharest, Romania, BSG Proceedings \textbf{12} (2005), Geometry Balkan Press, 252-258.
\bibitem{Viswanath-2000} Viswanath, D., \textit{Random Fibonacci sequences and the number} $1.13198824\ldots$, Math. Comput. \textbf{69}(231), 1131-1155. \\
\end{thebibliography}
\end{document}